\documentclass[11pt]{amsart}
\usepackage[a4paper, total={6in, 8in}]{geometry}
\usepackage{enumerate}
\usepackage{amsmath}
\usepackage{amssymb,latexsym}
\usepackage{amsthm}
\usepackage{color}
\usepackage{cancel}
\usepackage{graphicx}
\usepackage{cite}
\usepackage{changes}
\usepackage{lineno}
\usepackage{hyperref}
\usepackage{comment}
\usepackage{amscd}
\usepackage{cleveref}
\usepackage{mathrsfs}

\DeclareMathOperator{\C}{\mathcal{C}}

\DeclareMathOperator{\Aut}{Aut}
\DeclareMathOperator{\End}{End}
\DeclareMathOperator{\Gal}{Gal}

\DeclareMathOperator{\rk}{rk}

\newtheorem{theorem}{Theorem}[section]

\newtheorem{lemma}[theorem]{Lemma}
\newtheorem{corollary}[theorem]{Corollary}

\newtheorem{proposition}[theorem]{Proposition}

\newtheorem{example}[theorem]{Example}

\newtheorem{remark}[theorem]{Remark}

\newcommand{\F}{{\mathbb F}}

\newcommand{\SSS}{{\mathbb S}}
\newcommand{\NN}{{\mathbb N}}

\newcommand{\Tr}{{\mathrm Tr}}

\newcommand{\GL}{\hbox{{\rm GL}}}
\newcommand{\K}{{\mathbb K}}

\newcommand{\LL}{{\mathbb L}}

\newcommand{\PG}{\mathrm{PG}}
\newcommand{\N}{\mathrm{N}}
\newcommand{\gcrd}{\mathrm{gcrd}}
\newcommand{\lclm}{\mathrm{lclm}}

\newcommand{\Fix}{\mathrm{Fix}}

\newcommand{\EF}{\mathbb{E}_F}
\newcommand{\E}{\mathbb{E}}

\newcommand{\id}{\mathrm{id}}

\title{Iterative constructions of semifields 2-dimensional over a nucleus}
\author[L. Kölsch]{Lukas Kölsch}
\address{Lukas Kölsch,  \textnormal{Department of Mathematics and Statistics, University of South Florida,  140 7th Ave. South, St. Petersburg, FL 33701, USA \\}}
\email{lukas.koelsch.math@gmail.com}

\author[P. Santonastaso]{Paolo Santonastaso}
\address{Paolo Santonastaso, \textnormal{Department of Mathematics and Applications ``R. Caccioppoli'', University of Naples Federico II, Via Cintia, Monte Sant'Angelo, 80126 Naples, Italy \\}}
\email{paolo.santonastaso@unina.it}

\author[Y. Zhou]{Yue Zhou}
\address{Yue Zhou, \textnormal{College of Science, National University of Defense Technology, 410073 Changsha,
		China}}
\email{yue.zhou.ovgu@outlook.de}
\date{}

\begin{document}
	
	\begin{abstract}
		The construction of finite semifields dates back to Dickson's work in 1906. To the best of our knowledge, only two iterative constructions were previously known: those of Hiramine, Matsumoto and Oyama in 1987, and of Kantor and Williams in 2004. The latter yields, for suitable even orders, superpolynomially many pairwise non-isotopic semifields. In this paper, we use skew-polynomial rings to give a new iterative construction of finite semifields. This construction produces large families of pairwise non-isotopic semifields, each of which is $2$-dimensional over one of its nuclei. In particular, we show that this construction method yields semifields not contained in any family of semifields known so far.
	\end{abstract}
	
	\maketitle
    	{\textbf{Keywords}: semifield; iterative construction; skew polynomial; finite geometry}
	
	{\textbf{MSC2020}: 12K10; 51E15; 11T71; 16S36}
	
	
	\section{Introduction}
	Finite semifields are finite not-necessarily associative unital division algebras. Semifields have received a lot of attention in the last around 100 years, starting with the work of Dickson~\cite{dickson1906commutative} who constructed the first finite semifield that is not a field. Dickson's motivation to study semifields was purely of algebraic nature.  Subsequent inquiries into constructions and properties of semifields were often motivated by applications in both pure and applied mathematics, with semifields for instance exhibiting an intimate relation to non-Desarguesian projective planes in finite geometry~\cite{knuth1965finite} as well as optimal rank-metric codes~\cite{sheekeysurvey}. 
	
	The probably most fundamental question on finite semifields is about their ubiquity. Kantor~\cite{kantor_commutative_2003} asks if the number of non-isotopic semifields grows superpolynomially in the order of the semifield, providing an affirmative answer in the case of even order semifields~\cite{kantor2004symplectic}. In the odd order case, this question is still wide open. Indeed, the current best known bound is not even linear in the order of the semifield~\cite{gologlu2024counting}. This suggests that current constructions of odd order semifields are severely lacking. In particular, mimicking the approach of~\cite{kantor2004symplectic} that led to the superpolynomial count for even order semifields seems instructive. One key reason why Kantor and Williams~\cite{kantor2004symplectic} were able to achieve these counts was that their construction relied on an iterative approach, starting with a known semifield. 
	
	In this paper, we attempt a similar approach, realizing an iterative construction of semifields starting from a known family of semifields. We are however unable to achieve the superpolynomial count using our construction. However, our new construction still yields some very interesting results: Firstly, outside of the Kantor-Williams construction and the Hiramine-Matsumoto-Oyama (HMO) construction~\cite{hiramine1987extension}, our construction is the only iterative construction of semifields, opening the door to future refinements to attack Kantor's question on the number of odd order non-isotopic semifields. We expect that the techniques we develop here can yield crucial tools for a potential resolution of this question. Secondly, the semifields we construct have the peculiar property that they are always $2$-dimensional over a nucleus. While that is also probably a reason why the number of non-isotopic members in this family is not large enough to resolve Kantor's conjecture, the counts we do achieve show that the number of such semifields is still surprisingly large. Semifields two dimensional over a nucleus have a strong intrinsic structure (see the definition of nuclei later), for instance assuming some extra conditions on the other nuclei, semifields two dimensional over a nucleus have been completely classified~\cite{cardinali2006semifield}, with only a few families exhibiting these parameters. Our result shows that without the additional condition on the other  nuclei, large additional families are possible. Since nuclei lead for instance immediately to large autotopism groups (as well as collineation groups for the related projective planes), our construction thus yields many non-isotopic examples of highly symmetric semifields and semifield planes.
	
	After the preliminaries, \Cref{sec:skewpolynomialapproach}
	establishes our baseline construction and its nucleus calculations.
	The  iterated construction is then developed in
	Section~\ref{sec:iteration}. 
	In \Cref{sec:explicit-constructions}, we give explicit examples of our general construction, and \Cref{sec:linear-sets-isotopy}
	treats the isotopy question and establishes bounds on the number of non-isotopic semifields obtained from our construction.
	Finally, Section~\ref{sec:comparison-empt-hmo} compares these
	semifields with previously known constructions and proves that the
	family contains new isotopy classes.
	\section{Preliminaries}

	\subsection{Semifields} \label{sec:semifields}
	
	Let \(\F\) be a finite field. A finite-dimensional \(\F\)-algebra is a finite-dimensional \(\F\)-vector space \(\SSS\) endowed with an \(\F\)-bilinear map
	\[
	\star:\SSS\times \SSS\longrightarrow \SSS,
	\]
	which is not assumed to be associative. For \(a\in \SSS\), we denote by
	\[
	L_a:b\longmapsto a\star b,
	\qquad
	R_a:b\longmapsto b\star a
	\]
	the left and right multiplication maps. The algebra \((\SSS,+,\star)\) is called a \textbf{presemifield} if \(L_a\) and \(R_a\) are bijective for every \(a\neq 0\). If, in addition, there is an element \(1\in \SSS\) satisfying
	$
	1\star a=a\star 1=a$, for all $a\in\SSS$,
	then \(\SSS\) is called a \textbf{semifield}. Thus a finite semifield is a unital finite division algebra. Isotopy, introduced by Albert in \cite{albert1942non}, is the standard equivalence relation for finite semifields. Two presemifields \((\SSS,+,\star)\) and \((\SSS',+,\diamond)\) defined over the same prime field, are said to be \textbf{isotopic} if there exist additive bijections
	\[
	f,g,h:\SSS\longrightarrow \SSS'
	\]
	such that
	\[
	f(a\star b)=g(a)\diamond h(b),\]
	for all $a,b\in\SSS$. If one can take \(f=g=h\), then the two presemifields are isomorphic. By Kaplansky's trick, every finite presemifield is isotopic to a semifield. Therefore, throughout the rest of the paper, we may restrict our attention to semifields.

	Let \(\SSS=(\SSS,+,\star)\) be a semifield. Its left, middle, and right nuclei are defined by
	\[
	\NN_l(\SSS)
	:=
	\{a\in\SSS : a\star(b\star c)=(a\star b)\star c
	\text{ for all }b,c\in\SSS\},
	\]
	\[
	\NN_m(\SSS)
	:=
	\{b\in\SSS : a\star(b\star c)=(a\star b)\star c
	\text{ for all }a,c\in\SSS\},
	\]
	and
	\[
	\NN_r(\SSS)
	:=
	\{c\in\SSS : a\star(b\star c)=(a\star b)\star c
	\text{ for all }a,b\in\SSS\}.
	\]
	Their intersection
	\[
	\NN(\SSS):=\NN_l(\SSS)\cap \NN_m(\SSS)\cap \NN_r(\SSS)
	\]
	is the nucleus of \(\SSS\). The center is
	\[
	Z(\SSS)
	:=
	\{a\in\NN(\SSS): a\star b=b\star a
	\text{ for all }b\in\SSS\}.
	\]
	The nuclei and the center are finite fields.
	
	We shall use left multiplication to associate with \(\SSS\) its \textbf{spread set}. Put
	$
	\F:=Z(\SSS)$.
	Since every left multiplication map is \(\F\)-linear, we may consider
	\[
	\C(\SSS)
	:=
	\{L_a:a\in\SSS\}
	\subseteq \End_{\F}(\SSS).
	\]
	The map \(a\mapsto L_a\) is injective and $\F$-linear. Hence \(\C(\SSS)\) is an \(\F\)-linear subspace of \(\End_{\F}(\SSS)\), of dimension \(\dim_{\F}(\SSS)\), and every nonzero element of \(\C(\SSS)\) is nonsingular. Conversely, let \(V\) be a finite-dimensional vector space over a finite field \(\F\), and let
	$
	\C\subseteq \End_{\F}(V)$
	be an \(\F\)-subspace such that
	\[
	\dim_{\F}(\C)=\dim_{\F}(V),
	\qquad
	\det(T)\neq 0
	\quad\text{for every }0\neq T\in\C.
	\]
	After choosing an \(\F\)-linear isomorphism
	$
	\varphi:V\longrightarrow \C,
	$
	one obtains a presemifield multiplication on \(V\) by setting
	\[
	x\star_{\varphi} y:=\varphi(x)(y),
	\qquad x,y\in V.
	\]
	Different choices of \(\varphi\) give isotopic presemifields. 
	
	The spread set of a semifield can be also represented over its right nucleus. Indeed, \(\SSS\) is naturally a right vector space over \(\NN_r(\SSS)\), with scalar multiplication given by the semifield product. Moreover, for every \(a,b\in\SSS\) and every \(\alpha\in\NN_r(\SSS)\), one has
	\[
	L_a(b\star \alpha)
	=
	a\star(b\star\alpha)
	=
	(a\star b)\star\alpha
	=
	L_a(b)\star\alpha.
	\]
	Thus every \(L_a\) is \(\NN_r(\SSS)\)-linear, and hence
	\[
	\C(\SSS)
	\subseteq
	\End_{\NN_r(\SSS)}(\SSS)
	\subseteq
	\End_{\F}(\SSS).
	\]
	We shall use the following standard description of isotopy in terms of spread sets.
	
	\begin{proposition}
		\label{prop:isotopy-equivalence-right-nucleus}
		{\rm
			\cite[Theorem~7]{lavrauw2011finite}, \cite[Proposition~3]{sheekey2020new}}
		Let \(\SSS\) and \(\SSS'\) be semifields over the same right nucleus \(\E\), and let
		$
		\mathcal C(\SSS),\mathcal C(\SSS')\subseteq \End_{\E}(\SSS)
		$
		be their associated spread sets.
		Then \(\SSS\) and \(\SSS'\) are isotopic if and only if there exist
		\(A,B\in \GL_{\E}(\SSS)\) and \(\rho\in \Aut(\E)\) such that
		\[
		\mathcal C(\SSS')
		=
		A\,\mathcal C(\SSS)^{\rho}\,B.
		\]
		Here
		\[
		\mathcal C(\SSS)^{\rho}
		:=
		\{T^{\rho}:T\in \mathcal C(\SSS)\},
		\]
		where, after fixing an \(\E\)-basis of $\SSS$, \(\rho\) acts entrywise on matrices.
	\end{proposition}
	
	We now specialize to the case that is relevant for the present paper. Let \(\SSS\) be a semifield of order \(q^{2t}\) which is two-dimensional over its right nucleus. We identify
	\[
	\NN_r(\SSS)\cong \F_{q^t}
	\]
	and, as a right \(\F_{q^t}\)-vector space,
	\[
	\SSS\cong \F_{q^{2t}}.
	\]
	Therefore the spread set of \(\SSS\) may be regarded as a set
	\[
	\C(\SSS)
	\subseteq
	\End_{\F_{q^t}}(\F_{q^{2t}})
	\]
	of \(q^{2t}\) maps, containing the zero map and the identity map, with the property that every nonzero element is nonsingular. Every \(\F_{q^t}\)-linear endomorphism of \(\F_{q^{2t}}\) has a unique expression of the form
	\[
	T_{\alpha,\beta}:
	z\longmapsto \alpha z+\beta z^{q^t},
	\qquad
	\alpha,\beta\in\F_{q^{2t}}.
	\]
	
	The map $T_{\alpha,\beta}$ is nonsingular if and only if
	\[
	\mathrm{N}_{q^{2t}/q^t}(\alpha)
	\ne
	\mathrm{N}_{q^{2t}/q^t}(\beta).
	\]
	Indeed, suppose first that $\beta\ne 0$. A nonzero element $z$ lies
	in the kernel of $T_{\alpha,\beta}$ if and only if
	\[
	z^{q^t-1}=-\frac{\alpha}{\beta}.
	\]
	The image of the map $z\mapsto z^{q^t-1}$ on
	$\mathbb{F}_{q^{2t}}^\times$ is precisely the subgroup of elements
	of norm one over $\mathbb{F}_{q^t}$. Hence $T_{\alpha,\beta}$ has
	a nonzero kernel if and only if
	\[
	\mathrm{N}_{q^{2t}/q^t}\left(-\frac{\alpha}{\beta}\right)=1,
	\]
	or, equivalently,
	\[
	\mathrm{N}_{q^{2t}/q^t}(\alpha)
	=
	\mathrm{N}_{q^{2t}/q^t}(\beta).
	\]
	If $\beta=0$, then $T_{\alpha,0}$ is nonsingular if and only if
	$\alpha\ne 0$, which is again equivalent to the stated norm
	inequality, since $\mathrm{N}_{q^{2t}/q^t}(0)=0$. Consequently, the spread set of any semifield that is two-dimensional over its right nucleus can be written as
	\[
	S_{\mu,\zeta}
	:=
	\left\{
	T_{\mu(u),\zeta(u)}
	:
	u\in\F_{q^{2t}}
	\right\}
	\subseteq
	\End_{\F_{q^t}}(\F_{q^{2t}}),
	\]
	where \(\mu,\zeta:\F_{q^{2t}}\to\F_{q^{2t}}\) are additive maps such that
	$
	u\longmapsto (\mu(u),\zeta(u))$
	is injective, and each nonzero element of \(S_{\mu,\zeta}\) is nonsingular. Equivalently, for every $u\in\mathbb{F}_{q^{2t}}$ with
	$(\mu(u),\zeta(u))\ne(0,0)$, one has
	\[
	\mathrm{N}_{q^{2t}/q^t}\bigl(\mu(u)\bigr)
	\ne
	\mathrm{N}_{q^{2t}/q^t}\bigl(\zeta(u)\bigr).
	\]

	\subsection{Skew polynomial rings} \label{sec:introskew}
	
	Let $\K$ be a finite field and let $\LL / \K$ be a finite field extension with $[\LL : \K] = m$. 
	Let $\sigma$ be a generator of the Galois group $\Gal(\LL / \K)$, and consider the skew polynomial ring 
	\[
	R := \LL[X; \sigma].
	\]
	Thus, \(R\) consists of polynomials in the indeterminate \(X\) with coefficients in \(\LL\), written on the left, and multiplication is determined by the rule
	\[
	X\alpha=\sigma(\alpha)X,
	\]
	for all $\alpha\in \LL$, and extended distributively. The ring \(R\) is noncommutative as soon as \(\sigma\neq \mathrm{id}\). 
	Nevertheless, \(R\) retains many of the basic structural features of an ordinary polynomial ring. 
	In particular, every nonzero element \(f\in R\) has a well-defined degree, denoted by \(\deg(f)\), and for all nonzero \(f,g\in R\) one has
	\[
	\deg(fg)=\deg(f)+\deg(g),
	\qquad
	\deg(f+g)\leq \max\{\deg(f),\deg(g)\}.
	\]
	
	Moreover, \(R\) is both a left and a right Euclidean domain. In particular, given \(f,g\in R\) with \(g\neq 0\), there exist unique polynomials \(q,r\in R\) such that
	\[
	f=qg+r,
	\qquad\text{with either } r=0 \text{ or } \deg(r)<\deg(g).
	\]
	If \(r=0\), we say that \(g\) right-divides \(f\), and we write \(g\mid_r f\). Equivalently, \(f\) is a left multiple of \(g\). 
	As a consequence, every left ideal of \(R\) is principal. Hence, for any two nonzero elements \(f_1,f_2\in R\), one may define their (monic) greatest common right divisor and least common left multiple, denoted by
	\[
	\gcrd(f_1,f_2)
	\qquad\text{and}\qquad
	\lclm(f_1,f_2),
	\]
	respectively. If \(a=\gcrd(f_1,f_2)\) and \(b=\lclm(f_1,f_2)\), then
	\[
	Rf_1+Rf_2=Ra
	\qquad\text{and}\qquad
	Rf_1\cap Rf_2=Rb.
	\]
	It is well known that the center of \(R\) is
	\[
	Z(R)=\K[X^m].
	\]
	
	Now let \(F(Y)\in \K[Y]\) be a monic irreducible polynomial of degree \(s\ge 1\), with \(F(Y)\neq Y\). 
	Since \(F(X^m)\in Z(R)\), the ideal \(RF(X^m)\) is two-sided, and we may consider the quotient ring
	\[
	R_F:=\frac{R}{RF(X^m)}.
	\]
	Every element of \(R_F\) admits a unique representative of degree strictly smaller than \(ms\), so we identify each class \(\overline{a}\in R_F\) with its canonical representative of degree at most $ms-1$. Accordingly, when \(\overline{a}\neq 0\), we define
	$
	\deg(\overline{a}):=\deg(a),
	$
	where \(a\) is the canonical representative of \(\overline{a}\).
	
	The center of \(R_F\) will be denoted by \(\EF\). Since \(F(X^m)\) is central, one has
	\[
	\EF
	=
	Z(R_F)
	\cong
	\frac{\K[Y]}{\K[Y]F(Y)}.
	\]
	Because \(F(Y)\) is irreducible over \(\K\), it follows that \(\EF\) is a field extension of \(\K\) of degree \(s\), and hence
	\[
	\EF \cong \F_{|\K|^s}.
	\]
	
	Furthermore, the ideal \(RF(X^m)\) is maximal among the two-sided ideals of \(R\). Therefore \(R_F\) is a central simple algebra over \(\EF\), of dimension \(m^2\) over \(\EF\), and thus of dimension \(m^2s\) over \(\K\). By the Artin--Wedderburn theorem, one obtains an \(\EF\)-algebra isomorphism
	\begin{equation}\label{eq:artin}
		R_F
		\cong
		M_m(\EF).
	\end{equation}
	See, for instance, \cite[Theorem 1.2.19]{jacobson2009finite} and \cite[Theorem~20]{owen2023eigenspaces} for further details. 
	
	Fix an $\EF$-algebra isomorphism
	\[
	\varphi_F:R_F\longrightarrow M_m(\EF).
	\]
	For every $\bar a\in R_F$, we define
	\[
	\rk(\bar a):=
	\operatorname{rank}_{\EF}\bigl(\varphi_F(\bar a)\bigr).
	\] The following result expresses $\rk(\bar a)$ in terms of a
	greatest common right divisor in $R$.
	
	\begin{theorem}[see \textnormal{\cite[Proposition 7]{sheekey2020new}}]\label{th:rankpolynomial}
		Let \(\overline{a}=a+RF(X^m)\in R_F\) be nonzero. Then
			\[
			\rk(\bar a)
			=
			m-\frac{1}{s}\deg\bigl(\operatorname{gcrd}(a,F(X^m))\bigr).
			\]
	\end{theorem}

	\section{A lifting construction via skew polynomial rings}
	\label{sec:skewpolynomialapproach}
	
	In this section we describe the main lifting construction of the paper.  We start from a semifield of order \(q^{2t}\) which is two-dimensional over its right nucleus, and we use a skew-polynomial quotient to construct a semifield of order \(q^{2ts}\).  The idea is to embed the spread set of the initial semifield into the algebra
	\[
	R_G:=\frac{R}{RG(X^2)},
	\]
	where \(G\) is an irreducible polynomial over \(\F_{q^t}\). 
	
	Let \(n=2t\), with \(t\geq 1\), and put
	\[
	\K=\F_{q^t},
	\qquad
	\LL=\F_{q^n}=\F_{q^{2t}}.
	\]
	Let
	$\sigma:\LL\longrightarrow \LL$ be a 
	generator of \(\Gal(\LL/\K)\).  We consider the skew polynomial ring
	$
	R=\LL[X;\sigma]$.
	Since \([\LL:\K]=2\), the center of \(R\) is
	$
	Z(R)=\K[X^2]$.
	
	Let \(G(Y)\in \K[Y]\) be a monic irreducible polynomial of degree \(s\geq 1\), with \(G(Y)\neq Y\).  Since \(G(X^2)\in Z(R)\), the ideal \(RG(X^2)\) is two-sided, and the quotient
	\[
	R_G=\frac{R}{RG(X^2)}
	\]
	is a central simple algebra over
	\[
	Z(R_G)\cong \frac{\K[Y]}{\K[Y]G(Y)}\cong \F_{q^{ts}}.
	\]
	By \eqref{eq:artin}, we have an isomorphism of \(Z(R_G)\)-algebras
	\[
	R_G\cong M_2(\F_{q^{ts}}).
	\]
	We shall use the following rank-one criterion, which follows from \Cref{th:rankpolynomial} and \cite[Theorem~6]{sheekey2020new}.
	
	\begin{theorem} \label{th:normconditionsheekey} Let $\overline a = \alpha_0+\alpha_1X+\cdots+\alpha_sX^s+RG(X^2) \in R_G $ be nonzero. Then $ \rk(\overline a)\geq 1$. Moreover, if $ \rk(\overline a)=1,$ then $ \deg(\overline a)=s$ and \[ \frac{\N_{q^{2t}/q^t}(\alpha_0)} {\N_{q^{2t}/q^t}(\alpha_s)} = (-1)^sG(0). \] \end{theorem}

	\subsection{The lifted spread set}
	
	We now recall the starting point of the construction.  Let
	\[
	T(Y)=Y-1\in \K[Y].
	\]
	Then
	\[
	\frac{R}{RT(X^2)}
	=
	\frac{R}{R(X^2-1)}
	\cong M_2(\F_{q^t})
	\cong \End_{\F_{q^t}}(\F_{q^{2t}}).
	\]
	Under this identification, the class of \(\alpha+\beta X+RT(X^2)\) corresponds to the \(\F_{q^t}\)-linear map
	\[
	z\longmapsto \alpha z+\beta z^{q^t}.
	\]
	
	Let \(\mathcal S_{\mu,\zeta}\) be the spread set of a semifield of order \(q^{2t}\) which is two-dimensional over its right nucleus.  We write
	\[
	\mathcal S_{\mu,\zeta}
	=
	\left\{
	\mu(u)+\zeta(u)X+R(X^2-1)
	:
	u\in \F_{q^{2t}}
	\right\}
	\subseteq
	\frac{R}{R(X^2-1)},
	\]
	where
	\[
	\mu,\zeta:\F_{q^{2t}}\longrightarrow \F_{q^{2t}}
	\]
	are additive maps.  Since \(\mathcal S_{\mu,\zeta}\) is a spread set, the map
	\begin{equation}
		\label{eq:injectivesemifields}
		\Delta_{\mu,\zeta}:
		\F_{q^{2t}}\longrightarrow \F_{q^{2t}}^2,
		\qquad
		u\longmapsto (\mu(u),\zeta(u)),
	\end{equation}
	is injective, and every nonzero element of \(\mathcal S_{\mu,\zeta}\) is nonsingular.  Equivalently, for every \(u\in\F_{q^{2t}}\) with \(\zeta(u)\neq 0\), one has
	\begin{equation}
		\label{eq:norm-base-spread}
		\N_{q^{2t}/q^t}(\mu(u))
		\neq
		\N_{q^{2t}/q^t}(\zeta(u)).
	\end{equation}

	Let
	\[
	\eta\in \F_{q^{2t}}
	\]
	be such that
	\[
	\N_{\F_{q^{2t}}/\F_{q^t}}(\eta)=(-1)^sG(0).
	\]
	We define the lifted set associated with \(G\) and \(\eta\) as
	\begin{equation}
		\label{eq:extendingsemifield}
		\Phi_{G,\eta}\bigl(\mathcal S_{\mu,\zeta}\bigr)
		:=
		\left\{
		\mu(u_0)
		+
		\sum_{i=1}^{s-1}u_i X^i
		+
		\frac{\zeta(u_0)}{\eta}X^s+RG(X^2)
		:
		u_0,\ldots,u_{s-1}\in \F_{q^{2t}}
		\right\}
		\subseteq R_G.
	\end{equation}
	
	\begin{theorem}
		\label{th:liftingspread}
		If \(\mathcal S_{\mu,\zeta}\) is the spread set of a semifield of order \(q^{2t}\), then
		$
		\Phi_{G,\eta}\bigl(\mathcal S_{\mu,\zeta}\bigr)$
		is the spread set of a semifield of order \(q^{2ts}\).
	\end{theorem}
	
	\begin{proof}
		First, we prove that the lifted set has the correct size.  Since
		\[
		\mathcal S_{\mu,\zeta}
		=
		\left\{
		\mu(u)+\zeta(u)X+R(X^2-1):
		u\in \F_{q^{2t}}
		\right\}
		\]
		is a spread set, the map
		$
		\Delta_{\mu,\zeta}:
		u\longmapsto (\mu(u),\zeta(u))
		$
		is injective.  Hence the parameters
		$
		u_0,u_1,\ldots,u_{s-1}\in \F_{q^{2t}}
		$
		determine uniquely an element of
		$
		\Phi_{G,\eta}\bigl(\mathcal S_{\mu,\zeta}\bigr)$.
		Therefore
		\[
		\left|
		\Phi_{G,\eta}\bigl(\mathcal S_{\mu,\zeta}\bigr)
		\right|
		=
		q^{2ts}.
		\]
		Moreover, the set is additive, because \(\mu\) and \(\zeta\) are additive and the coefficients \(u_1,\ldots,u_{s-1}\) run over $\F_{q^n}$. Since the original spread set contains the identity map, there exists \(e\in \F_{q^{2t}}\) such that
		$
		\mu(e)=1$ and $\zeta(e)=0$.
		Hence
		\[
		1+RG(X^2)\in \Phi_{G,\eta}\bigl(\mathcal S_{\mu,\zeta}\bigr).
		\]
		It remains to prove that every nonzero element of
		$
		\Phi_{G,\eta}\bigl(\mathcal S_{\mu,\zeta}\bigr)
		$
		is nonsingular in
		$
		R_G\cong M_2(\F_{q^{ts}})$.
		Suppose, by contradiction, that there exists a nonzero singular element
		\[
		\overline{a}
		=
		\mu(u_0)
		+
		\sum_{i=1}^{s-1}u_i X^i
		+
		\frac{\zeta(u_0)}{\eta}X^s
		+
		RG(X^2)
		\]
		in the lifted set.  Since \(\overline{a}\neq 0\), we have
		$
		\rk(\overline{a})\geq 1.
		$
		As \(\overline{a}\) is singular in \(M_2(\F_{q^{ts}})\), we must have
		\[
		\rk(\overline{a})=1.
		\]
		Therefore, by \Cref{th:normconditionsheekey}, the coefficient of \(X^s\) is nonzero, so
		$
		\zeta(u_0)\neq 0,$
		and
		\[
		\frac{
			\N_{\F_{q^{2t}}/\F_{q^t}}(\mu(u_0))
		}{
			\N_{\F_{q^{2t}}/\F_{q^t}}\!\left(\zeta(u_0)/\eta\right)
		}
		=
		(-1)^sG(0).
		\]
		Since
		\[
		\N_{\F_{q^{2t}}/\F_{q^t}}\!\left(\frac{\zeta(u_0)}{\eta}\right)
		=
		\frac{
			\N_{\F_{q^{2t}}/\F_{q^t}}(\zeta(u_0))
		}{
			\N_{\F_{q^{2t}}/\F_{q^t}}(\eta)
		},
		\]
		we get
		\[
		\frac{
			\N_{\F_{q^{2t}}/\F_{q^t}}(\mu(u_0))
		}{
			\N_{\F_{q^{2t}}/\F_{q^t}}(\zeta(u_0))
		}
		\,
		\N_{\F_{q^{2t}}/\F_{q^t}}(\eta)
		=
		(-1)^sG(0).
		\]
		Using
		$
		\N_{\F_{q^{2t}}/\F_{q^t}}(\eta)=(-1)^sG(0),
		$
		we obtain
		\[
		\N_{\F_{q^{2t}}/\F_{q^t}}(\mu(u_0))
		=
		\N_{\F_{q^{2t}}/\F_{q^t}}(\zeta(u_0)).
		\]
		But \(\zeta(u_0)\neq 0\), and this contradicts the nonsingularity of the original spread set
		$
		\mathcal S_{\mu,\zeta}$. Therefore, every nonzero element of
		$
		\Phi_{G,\eta}\bigl(\mathcal S_{\mu,\zeta}\bigr)
		$
		is nonsingular.  Thus
		$
		\Phi_{G,\eta}\bigl(\mathcal S_{\mu,\zeta}\bigr)
		$
		is an additive spread set of size \(q^{2ts}\), and hence it defines a semifield of order \(q^{2ts}\).
	\end{proof}
	
	\begin{remark}\label{rem:lifting-presemifields-quasifields}
			The construction also applies to presemifield spread sets in
			$\operatorname{End}_{\mathbb{F}_{q^t}}
			(\mathbb{F}_{q^{2t}})$.
			Indeed, the assumption that the initial spread set contains the
			identity is used only to ensure that the lifted set also contains
			the identity. Without this assumption, the same proof shows that
			the lifted set is an additive spread set of cardinality $q^{2ts}$,
			and hence defines a presemifield of order $q^{2ts}$.
			
			More generally, the additivity of $\mu$ and $\zeta$ can be dropped
			provided that the difference of any two distinct elements of
			$S_{\mu,\zeta}$ is nonsingular. Explicitly, suppose that
			$\mu,\zeta:\mathbb{F}_{q^{2t}}\to\mathbb{F}_{q^{2t}}$ satisfy
			$\mu(0)=\zeta(0)=0$ and
			\[
			\mathrm{N}_{q^{2t}/q^t}\bigl(\mu(u)-\mu(v)\bigr)
			\ne
			\mathrm{N}_{q^{2t}/q^t}\bigl(\zeta(u)-\zeta(v)\bigr)
			\qquad
			\text{for all }u\ne v.
			\]
			Then the same formula defining $\Phi_{G,\eta}(S_{\mu,\zeta})$
			produces a set of cardinality $q^{2ts}$ containing zero, in which
			the difference of any two distinct elements is nonsingular. Thus the construction extends to quasifield spread sets in
			$M_2(\mathbb{F}_{q^t})$, producing spread sets in
			$M_2(\mathbb{F}_{q^{ts}})$ and, after normalization if necessary,
			quasifields of order $q^{2ts}$.
	\end{remark}
	
	\begin{remark}
		Notice that, when \(s\) is odd, the choice of the element \(\eta\) satisfying
		the prescribed norm condition does not affect the isotopy class of the
		resulting semifield. Indeed, let \(\eta,\eta'\in\F_{q^{2t}}^\times\) be such that
		\[
		\N_{\F_{q^{2t}}/\F_{q^t}}(\eta)
		=
		\N_{\F_{q^{2t}}/\F_{q^t}}(\eta')
		=
		(-1)^sG(0).
		\]
		Then \(d=\eta/\eta'\) has norm one over \(\F_{q^t}\), and hence, by Hilbert's
		Theorem~90, there exists \(c\in \F_{q^{2t}}^\times\) such that
		\[
		d=c^{-1}\sigma(c).
		\]
		Conjugation by \(c\) in \(R_G\) fixes the even powers of \(X\), while it
		multiplies the coefficients of the odd powers by \(d\). Since \(s\) is odd,
		the coefficient of \(X^s\) is therefore transformed as
		\[
		\frac{\zeta(u_0)}{\eta}
		\longmapsto
		d\,\frac{\zeta(u_0)}{\eta}
		=
		\frac{\zeta(u_0)}{\eta'}.
		\]
		On the other hand, the coefficients of \(X,\ldots,X^{s-1}\) range
		independently over \(\F_{q^{2t}}\), so multiplying the odd-indexed ones by the
		nonzero scalar \(d\) does not change their set of possible values. It follows
		that
		\[
		c^{-1}\Phi_{G,\eta}\bigl(\mathcal S_{\mu,\zeta}\bigr)c
		=
		\Phi_{G,\eta'}\bigl(\mathcal S_{\mu,\zeta}\bigr).
		\]
		Thus different admissible choices of \(\eta\) yield equivalent spread sets,
		and hence isotopic semifields.
	\end{remark}

	\subsection{The left and middle nuclei}
	
	In the following, we compute the left and middle nuclei of the semifields obtained from the lifted spread set
	$
	\Phi_{G,\eta}\bigl(\mathcal S_{\mu,\zeta}\bigr)$.
	These computations can be performed directly on the spread set by using the standard idealiser description of nuclei.
	
	\begin{theorem} [see \textnormal{\cite[Theorem 1]{marino2012nuclei} and \cite[Proposition 5]{sheekey2020new}}]
		\label{th:thmaps}
		Let \(\SSS\) be a semifield and let
		$
		\mathcal C(\SSS)\subseteq \End_{\F}(\SSS)
		$
		be its spread set.  Then
		\[
		\NN_l(\SSS)
		\simeq
		\left\{
		\varphi\in \End_{\F}(\SSS):
		\varphi\circ \psi\in \mathcal C(\SSS)
		\text{ for every }\psi\in \mathcal C(\SSS)
		\right\},
		\]
		and
		\[
		\NN_m(\SSS)
		\simeq
		\left\{
		\varphi\in \End_{\F}(\SSS):
		\psi\circ \varphi\in \mathcal C(\SSS)
		\text{ for every }\psi\in \mathcal C(\SSS)
		\right\}.
		\]
	\end{theorem}
	
	Let
	\[
	\mathcal C_{G,\eta}
	:=
	\Phi_{G,\eta}\bigl(\mathcal S_{\mu,\zeta}\bigr).
	\]
	Thus
	\[
	\mathcal C_{G,\eta}
	=
	\left\{
	\mu(u_0)
	+
	\sum_{i=1}^{s-1}u_i X^i
	+
	\frac{\zeta(u_0)}{\eta}X^s
	+
	RG(X^2)
	:
	u_0,\ldots,u_{s-1}\in \F_{q^{2t}}
	\right\}.
	\]

	\begin{theorem}
		\label{th:left-middle-nuclei-lift}
		Assume \(s\geq 2\), and let \(\SSS_{G,\eta}\) be the semifield associated with the spread set \(\mathcal C_{G,\eta}\).  Then
		\[
		\NN_l(\SSS_{G,\eta})
		\simeq
		\mathcal N_l(\mu,\zeta),
		\]
		where
		\[
		\mathcal N_l(\mu,\zeta)
		=
		\left\{
		\alpha\in \F_{q^{2t}}
		:
		\forall u\in \F_{q^{2t}},
		\ \exists v\in \F_{q^{2t}}
		\text{ such that }
		(\mu(v),\zeta(v))
		=
		(\alpha\mu(u),\alpha\zeta(u))
		\right\}.
		\]
		Moreover,
		\[
		\NN_m(\SSS_{G,\eta})
		\simeq
		\mathcal N_m^{(s)}(\mu,\zeta),
		\]
		where
		\[
		\mathcal N_m^{(s)}(\mu,\zeta)
		=
		\left\{
		\alpha\in \F_{q^{2t}}
		:
		\forall u\in \F_{q^{2t}},
		\ \exists v\in \F_{q^{2t}}
		\text{ such that }
		(\mu(v),\zeta(v))
		=
		(\mu(u)\alpha,\zeta(u)\sigma^s(\alpha))
		\right\}.
		\]
	\end{theorem}

	\begin{proof}
		We use \Cref{th:thmaps} over the prime field.  Since
		$
		1\in \mathcal C_{G,\eta},$
		every element of the left idealiser
		\[
		\{h\in R_G:h\mathcal C_{G,\eta}\subseteq \mathcal C_{G,\eta}\}
		\]
		belongs to \(\mathcal C_{G,\eta}\).  Hence such an element has a representative of the form
		\[
		\overline{h}=h_0+h_1X+\cdots+h_sX^s+RG(X^2).
		\]
		
		We first prove that \(\overline{h}\) must be a scalar.  Since \(s\geq 2\), we have
		\[
		X+RG(X^2)\in \mathcal C_{G,\eta}.
		\]
		Thus
		$
		\overline{h}X\in \mathcal C_{G,\eta}$. Then
		\[
		\overline{h}X=h_0X+h_1X^2+\cdots+h_sX^{s+1}\in \mathcal C_{G,\eta}.
		\]
		Since \(s+1<2s\), the coefficient of \(X^{s+1}\) is not affected by reduction modulo \(G(X^2)\).  But every element of \(\mathcal C_{G,\eta}\) has degree at most \(s\). Therefore
		$
		h_s=0$.
		Hence
		\[
		\deg(\overline{h})\leq s-1.
		\]
		Moreover,
		$
		X^{s-1}+RG(X^2)\in \mathcal C_{G,\eta}.
		$
		Therefore
		\[
		\overline{h}X^{s-1}\in \mathcal C_{G,\eta}.
		\]
		The terms
		\[
		h_i X^{i+s-1},
		\qquad
		2\leq i\leq s-1,
		\]
		have degrees strictly larger than \(s\) and strictly smaller than \(2s\).  Again, they cannot occur in an element of \(\mathcal C_{G,\eta}\).  Thus
		$
		h_i=0$, for every $2\leq i\leq s-1$.
		Hence
		\[
		\overline{h}=h_0+h_1X+RG(X^2).
		\]
		
		In the next step we show that \(h_1=0\).  Since the initial semifield is two-dimensional over its right nucleus, it is not a field.  Hence \(\zeta\) is not identically zero.  Choose \(u\in \F_{q^{2t}}\) such that
		\[
		\zeta(u)\neq 0.
		\]
		Then
		\[
		\overline{a}
		=
		\mu(u)+\frac{\zeta(u)}{\eta}X^s+RG(X^2)
		\]
		belongs to \(\mathcal C_{G,\eta}\).  Since
		$
		h\mathcal C_{G,\eta}\subseteq \mathcal C_{G,\eta},$
		we have
		$
		\overline{h}\overline{a}\in \mathcal C_{G,\eta}.
		$
		The coefficient of \(X^{s+1}\) in
		\[
		(h_0+h_1X)
		\left(
		\mu(u)+\frac{\zeta(u)}{\eta}X^s
		\right)
		\]
		is
		\[
		h_1\sigma\!\left(\frac{\zeta(u)}{\eta}\right).
		\]
		Since \(\zeta(u)\neq 0\), this coefficient vanishes if and only if
		\[
		h_1=0.
		\]
		Therefore
		\[
		\overline{h}=h_0+RG(X^2),
		\]
		with $h_0 \in \F_{q^{2t}}$.
		
		Thus the left nucleus is represented by those scalars \(\alpha\in \F_{q^{2t}}\) such that
		\[
		\alpha \overline{a}\in \mathcal C_{G,\eta}
		\]
		for every
		$
		\overline{a}\in \mathcal C_{G,\eta}$.
		Let
		\[
		\overline{a}
		=
		\mu(u)
		+
		\sum_{i=1}^{s-1}u_i X^i
		+
		\frac{\zeta(u)}{\eta}X^s+RG(X^2) \in \mathcal C_{G,\eta}.
		\]
		Then
		\[
		\alpha \overline{a}
		=
		\alpha\mu(u)
		+
		\sum_{i=1}^{s-1}\alpha u_i X^i
		+
		\frac{\alpha\zeta(u)}{\eta}X^s+RG(X^2).
		\]
		Since the coefficients of \(X,\ldots,X^{s-1}\) are free in the definition of \(\mathcal C_{G,\eta}\), we have
		$
		\alpha \overline{a}\in \mathcal C_{G,\eta}
		$
		for every \(\overline{a}\in \mathcal C_{G,\eta}\) if and only if, for every \(u\in \F_{q^{2t}}\), there exists \(v\in \F_{q^{2t}}\) such that
		\[
		\mu(v)=\alpha\mu(u)
		\]
		and
		\[
		\zeta(v)=\alpha\zeta(u).
		\]
		This gives the asserted description of the left nucleus. The proof for the middle nucleus is similar.
	\end{proof}

	\subsection{The right nucleus}
	
	We now determine the right nucleus of the semifield associated
	with the lifted spread set $\C_{G,\eta}$. Throughout this subsection,
	we assume that $s\geq 2$.
	
	Recall that
	\[
	E_G:=Z(R_G)
	\cong \frac{\K[Y]}{\K[Y]G(Y)}
	\cong \mathbb F_{q^{ts}}.
	\]
	Put $V=E_G^2$, and consider an $E_G$-algebra isomorphism
	\[
	\varphi_G\colon R_G
	\longrightarrow M_2(E_G)
	\cong \operatorname{End}_{E_G}(V).
	\]
	In particular,
	$
	\varphi_G(\lambda)=\lambda I_V$, for every $\lambda\in E_G.$ By Theorem \ref{th:liftingspread}, the set $\varphi_G(\C_{G,\eta})$ is an
	additive spread set containing $I_V$. As in the previous
	subsection, let
	\[
	\SSS_{G,\eta}=(V,+,\star_{G,\eta})
	\]
	be a semifield with spread set
	\[
	\C(\SSS_{G,\eta})=\varphi_G(\C_{G,\eta}).
	\] Denote its multiplicative identity by $e$. Following the
	notation of Subsection~\ref{sec:semifields}, write
	\[
	L_x(y)=x\star_{G,\eta}y,
	\qquad
	R_x(y)=y\star_{G,\eta}x
	\qquad (x,y\in V).
	\]
	
	\begin{theorem}\label{th:right-nucleus-lift}
		The right nucleus of $\SSS_{G,\eta}$ is
		\[
		\NN_r(\SSS_{G,\eta})
		=
		\{\lambda e:\lambda\in E_G\}
		\cong E_G
		\cong \mathbb F_{q^{ts}}.
		\]
		Consequently, $\SSS_{G,\eta}$ is two-dimensional over its
		right nucleus.
	\end{theorem}
	
	\begin{proof}
		Let $\mathbb F_p$ be the prime subfield of $\K$. Consider
		the centraliser of $\varphi_G(\C_{G,\eta})$ in
		$\operatorname{End}_{\mathbb F_p}(V)$:
		\[
		\mathcal Z
		:=
		\left\{
		D\in\operatorname{End}_{\mathbb F_p}(V):
		D\circ\varphi_G(\bar a)
		=
		\varphi_G(\bar a)\circ D
		\text{ for every }\bar a\in \C_{G,\eta}
		\right\}.
		\]
		The endomorphisms in $\mathcal Z$ are initially assumed
		to be only $\mathbb F_p$-linear. By the centraliser description of the right nucleus,
		applied over $\mathbb F_p$,
		\[
		\mathcal Z
		=
		\{R_c:c\in \NN_r(\SSS_{G,\eta})\} \cong \NN_r(\SSS_{G,\eta});
		\]
		see \cite[Proposition~5(iii)]{sheekey2020new}.
		
		We now prove that
		\[
		\mathcal Z=\{\lambda I_V:\lambda\in E_G\}.
		\]
		Since every element of $\varphi_G(\C_{G,\eta})$ is
		$E_G$-linear, every scalar endomorphism $\lambda I_V$,
		with $\lambda\in E_G$, belongs to $\mathcal Z$.
		
		For the reverse inclusion, let $D\in\mathcal Z$.
		For $f\in R$, write
		$
		\bar f:=f+RG(X^2)\in R_G$.
		Since $s\geq 2$ and $\mu(0)=\zeta(0)=0$, the definition
		of $\C_{G,\eta}$ gives
		\[
		\overline{aX}=\bar a\,\bar X\in \mathcal{C}_{G,\eta}
		\qquad\text{for every }a\in \LL.
		\]
		In particular, $\bar X$ is a nonzero element of
		$\C_{G,\eta}$, so $\varphi_G(\bar X)$ is invertible.
		As $D$ commutes with $\varphi_G(\bar X)$, it also
		commutes with $\varphi_G(\bar X)^{-1}$. For every $a\in \LL$, we have
		\[
		\varphi_G(\bar a)
		=
		\varphi_G(\overline{aX})
		\circ\varphi_G(\bar X)^{-1}.
		\]
		Both factors on the right commute with $D$.
		Consequently,
		\begin{equation}\label{eq:right_nucleus_D_circ}
			D\circ\varphi_G(\bar a)
			=
			\varphi_G(\bar a)\circ D
			\qquad\text{for every }a\in \LL.	
		\end{equation}
		
		Every element $\bar f\in R_G$ has a unique representative
		of the form
		\[
		f=\sum_{j=0}^{2s-1}a_jX^j,
		\qquad a_j\in \LL,
		\]
		and hence
		\[
		\varphi_G(\bar f)
		=
		\sum_{j=0}^{2s-1}
		\varphi_G(\bar a_j)
		\circ\varphi_G(\bar X)^j.
		\]
		Since $D$ is additive and commutes with
		$\varphi_G(\bar X)$ and with every $\varphi_G(\bar a_j)$,
		it follows from \eqref{eq:right_nucleus_D_circ} that
		\[
		D\circ\varphi_G(\bar f)
		=
		\varphi_G(\bar f)\circ D
		\qquad\text{for every }\bar f\in R_G.
		\]
		In particular, $D$ commutes with
		$\varphi_G(\lambda)=\lambda I_V$ for every $\lambda\in E_G$.
		Therefore
		\[
		D(\lambda v)=\lambda D(v)
		\qquad (\lambda\in E_G,\ v\in V).
		\]
		Thus $D$ is $E_G$-linear. 
		
		We may now use the surjectivity of $\varphi_G$ to write
		$
		D=\varphi_G(\bar d)$
		for a unique $\bar d\in R_G$. The commutation relations
		above imply
		\[
		\varphi_G(\bar d\,\bar f)
		=
		\varphi_G(\bar f\,\bar d)
		\qquad\text{for every }\bar f\in R_G.
		\]
		By injectivity of $\varphi_G$,
		\[
		\bar d\,\bar f=\bar f\,\bar d
		\qquad\text{for every }\bar f\in R_G.
		\]
		Hence $\bar d\in Z(R_G)=E_G$, and therefore
		$D=\lambda I_V$ for some $\lambda\in E_G$.
		This proves
		\[
		\mathcal Z=\{\lambda I_V:\lambda\in E_G\}\cong E_G.
		\]
		Since $\dim_{E_G}(V)=2$, the semifield $\SSS_{G,\eta}$
		is two-dimensional over its right nucleus.
	\end{proof}

	\section{Iterative constructions}
	\label{sec:iteration}
	In this section we show that the lifting construction developed above can be iterated. The main point is that, after one application of the construction, the resulting spread set can again be represented in a form suitable for applying the same procedure. To make this explicit, we use a semilinear transformation model for the quotient algebra arising from the skew polynomial construction.

	\subsection{Irreducible semilinear maps}
	\label{subsec:irreducible-semilinear-maps}
	
	We now describe a semilinear transformations framework corresponding to the skew-polynomial
	algebra used in Section~\ref{sec:skewpolynomialapproach}. We keep the notation
	\[
	\K=\F_{q^t},
	\qquad
	\LL=\F_{q^{2t}},
	\]
	and we denote by
	\[
	\sigma:\LL\longrightarrow \LL,
	\qquad
	a\longmapsto a^{q^t},
	\]
	the non-trivial element of $\Gal(\LL/\K)$. Thus $\sigma$ has order $2$ and, for
	the skew polynomial ring
	$
	R=\LL[X;\sigma],
	$
	one has
	$
	Z(R)=\K[X^2]$.
	
	Let $V$ be a vector space over $\LL$. A map
	$
	T:V\longrightarrow V$
	is called $\sigma$-\textbf{semilinear} if it is additive and satisfies
	$
	T(\alpha v)=\sigma(\alpha)T(v),$
	for every $\alpha\in\LL$ and $v\in V$. The study of semilinear maps goes back to
	Jacobson~\cite{jacobson1937pseudo}. Of particular relevance in the theory of
	finite semifields are the irreducible semilinear maps; see, for instance,
	\cite{kantor2008semifields,lavrauw2013semifields,Dempwolff2010onirreducible}. A $\sigma$-semilinear map $T$ is said to be \textbf{irreducible} if the only
	$T$-invariant $\LL$-subspaces of $V$ are $\{0\}$ and $V$. In
	\cite{Dempwolff2010onirreducible}, Dempwolff gives a characterization of
	irreducible semilinear transformations and an explicit description of them. We
	now recall the special case needed in this paper.
	
	Let $s\geq 3$ be odd, and let
	\[
	G(Y)\in\K[Y]
	\]
	be a monic irreducible polynomial of degree $s$, with $G(Y)\neq Y$. Let $u$ be
	a root of $G(Y)$ and put
	\[
	\E=\K(u)\cong \frac{\K[Y]}{\K[Y]G(Y)}\cong \F_{q^{ts}}.
	\]
	Since $\gcd(s,2)=1$, the fields $\LL$ and $\E$ are linearly disjoint over
	$\K$. Hence their compositum is
	\[
	V:=\LL\E=\F_{q^{2ts}}.
	\]
	In particular, $V$ is an $s$-dimensional vector space over $\LL$ and a
	two-dimensional vector space over $\E$.
	
	Let
	\[
	\gamma:V\longrightarrow V,
	\qquad
	x\longmapsto x^{q^{st}}.
	\]
	Then $\gamma$ is the non-trivial element of $\Gal(V/\E)$. Moreover, since $s$
	is odd, the restriction of $\gamma$ to $\LL=\F_{q^{2t}}$ is precisely $\sigma$.
	Indeed,
	$
	st\equiv t \pmod {2t}$.
	Thus
	\[
	\gamma_{|\LL}=\sigma.
	\]
	Since the norm map
	$
	\N_{V/\E}:V^\times\longrightarrow \E^\times$
	is surjective, we may choose $w\in V^\times$ such that
	\[
	\N_{V/\E}(w)=w\gamma(w)=u.
	\]
	Define
	\[
	T_w:V\longrightarrow V,
	\qquad
	x\longmapsto wx^{\gamma}.
	\]
	Then $T_w$ is $\sigma$-semilinear with respect to the $\LL$-vector space
	structure of $V$. Indeed, for every $\alpha \in\LL$ and $x\in V$, we have
	\[
	T_w(\alpha x)
	=
	w\gamma(\alpha)x^{\gamma}
	=
	\sigma(\alpha)T_w(x).
	\]
	Moreover, $T_w$ is $\E$-linear, because $\gamma$ fixes $\E$ pointwise. Finally,
	\[
	T_w^2(x)
	=
	w\gamma(w)x^{\gamma^2}
	=
	ux
	\]
	for every $x\in V$. Hence
	\[
	T_w^2=u\,\id_V.
	\]
	Since $G(u)=0$, we get
	\[
	G(T_w^2)=0.
	\]
	
	\begin{proposition}
		\label{prop:Tw-irreducible}
		The map $T_w$ is an irreducible $\sigma$-semilinear map on the
		$\LL$-vector space $V$.
	\end{proposition}
	
	\begin{proof}
		Let $U$ be a non-zero $\LL$-subspace of $V$ such that
		\[
		T_w(U)\subseteq U.
		\]
		Then
		\[
		T_w^2(U)\subseteq U.
		\]
		Since $T_w^2$ is multiplication by $u$, the subspace $U$ is stable under
		multiplication by $u$. It is also stable under multiplication by every element
		of $\LL$, because $U$ is an $\LL$-subspace. Therefore $U$ is stable under
		multiplication by the field generated by $\LL$ and $u$, namely
		\[
		\LL(u)=\LL\E=V.
		\]
		Thus $U$ is a non-zero $V$-subspace of the one-dimensional $V$-vector space
		$V$. Hence $U=V$. This proves that $T_w$ has no non-trivial invariant
		$\LL$-subspaces, and so $T_w$ is irreducible.
	\end{proof}
	
	\begin{remark}
		The irreducibility of $T_w$ can also be deduced from Dempwolff's description of
		irreducible semilinear transformations. Indeed, in the notation of
		\cite[Section~2.8]{Dempwolff2010onirreducible}, we are considering the case
		$m=2$ and $\gcd(2,s)=1$. The condition
		\[
		\K\bigl(\N_{V/\E}(w)\bigr)=\E
		\]
		implies that $T_w$ is irreducible. In our situation
		$\N_{V/\E}(w)=u$, and $G(Y)$ is the minimal polynomial of $u$ over $\K$, so
		$\K(u)=\E$.
	\end{remark}
	
	We now connect this semilinear model with the skew-polynomial quotient used in
	the lifting construction. Recall that
	\[
	R_G=\frac{R}{RG(X^2)}.
	\]
	We identify the center of $R_G$ with $\E$ through the isomorphism
	\[
	\frac{\K[Y]}{\K[Y]G(Y)}
	\longrightarrow
	\E,
	\qquad
	Y+\K[Y]G(Y)\longmapsto u.
	\]
	
	Since $\gamma$ has order $2$, we have
	\[
	T_w^2(x)=ux,
	\]
	and hence,  for $i\geq 0$,
	\[
	T_w^{2i}(x)=u^i x,
	\qquad
	T_w^{2i+1}(x)=u^i wx^{\gamma},
	\]

	\begin{theorem} [see \textnormal{\cite{santonastaso2026autotopism}}]
		\label{th:SZ-specialized-RG}
		Let $G(Y)\in\K[Y]$ be a monic irreducible polynomial of odd degree $s$, with
		$G(Y)\neq Y$. Let $u\in\E$ be a root of $G(Y)$, and let $w\in V^\times$ be such
		that
		$
		\N_{V/\E}(w)=u$.
		Then the map
		\begin{equation} \label{eq:fromskewtosemilinear}
			\begin{array}{rccc}
				\Psi_w: & R_G & \longrightarrow & \End_{\E}(V)\\[4pt]
				& \displaystyle
				\sum_{i=0}^{2s-1}\alpha_iX^i+RG(X^2)
				& \longmapsto &
				\displaystyle \sum_{i=0}^{2s-1}\alpha_i  T_{w}^i
			\end{array}
		\end{equation}
		is an isomorphism of $\E$-algebras. 
	\end{theorem}
	
	We note that if
	\[
	\overline a
	=
	\sum_{i=0}^{s-1}\alpha_iX^{2i}
	+
	\sum_{i=0}^{s-1}\beta_iX^{2i+1}
	+
	RG(X^2),
	\]
	with $\alpha_i,\beta_i\in\LL$, then
	\[
	\Psi_w(\overline a)(x)
	=
	\sum_{i=0}^{s-1}\alpha_i u^i x
	+
	\sum_{i=0}^{s-1}\beta_i u^i wx^{\gamma},
	\]
	for every $x\in V$.  Under the isomorphism $\Psi_w$, the lifted spread set
	$
	\Phi_{G,\eta}\bigl(\mathcal S_{\mu,\zeta}\bigr)
	\subseteq R_G
	$
	given in \eqref{eq:extendingsemifield}
	can be described explicitly as a set of $\E$-linear endomorphisms of $V$.
	
	Let
	\[
	s=2r+1.
	\]
	For
	\[
	\overline a
	=
	\mu(u_0)
	+
	\sum_{i=1}^{s-1}u_iX^i
	+
	\frac{\zeta(u_0)}{\eta}X^s
	+
	RG(X^2),
	\]
	with $u_0,\ldots,u_{s-1}\in\LL$, we have
	\[
	\begin{array}{rl}
		\Psi_w(\overline a)(x)
		& =
		\mu(u_0)x
		+
		\sum\limits_{i=1}^{s-1}u_iT_w^i(x)
		+
		\frac{\zeta(u_0)}{\eta}T_w^s(x) \\
		& =
		\left(
		\mu(u_0)
		+
		\sum\limits_{j=1}^{r}u_{2j}u^j
		\right)x
		+
		\left(
		\sum\limits_{j=0}^{r-1}u_{2j+1}u^j
		+
		\frac{\zeta(u_0)}{\eta}u^r
		\right)wx^\gamma.
	\end{array}
	\]
	Therefore
	\[
	\Psi_w\left(
	\Phi_{G,\eta}\bigl(\mathcal S_{\mu,\zeta}\bigr)
	\right)
	=
	\left\{
	A(u_0,\ldots,u_{s-1})x
	+
	B(u_0,\ldots,u_{s-1}) wx^{\gamma}
	:
	u_0,\ldots,u_{s-1}\in\LL
	\right\},
	\]
	where
	\[
	A(u_0,\ldots,u_{s-1})
	=
	\mu(u_0)
	+
	\sum_{j=1}^{r}u_{2j}u^j
	\]
	and
	\[
	B(u_0,\ldots,u_{s-1})
	=
	\sum_{j=0}^{r-1}u_{2j+1}u^j
	+
	\frac{\zeta(u_0)}{\eta}u^r.
	\]
	Thus the lifted spread set is represented, in the semilinear model, by the maps
	\[
	x\longmapsto A x+Bwx^\gamma,
	\]
	where $A$ and $B$ are constrained by the above parametrization. Now, Theorem~\ref{th:liftingspread} shows that, whenever
	$\mathcal S_{\mu,\zeta}$ is the spread set of a semifield, the set
	\[
	\Psi_w\left(
	\Phi_{G,\eta}\bigl(\mathcal S_{\mu,\zeta}\bigr)
	\right)
	\]
	is again the spread set of a semifield. Furthermore, by
	\Cref{th:right-nucleus-lift}, the resulting semifield is two-dimensional over
	its right nucleus. Hence the lifting construction can be iterated, and this is
	the key idea developed in the next subsection.

	\subsection{Iterating the lifting construction}
	\label{subsec:iterating-lifting-construction}
	
	We now explain how the previous construction can be iterated. The main point is
	that, after one lifting step, the resulting semifield is again
	two-dimensional over its right nucleus. Hence it can be used as input for the
	same construction.
	
	Let
	\[
	t_0=t,
	\qquad
	\K_0=\F_{q^{t_0}},
	\qquad
	\LL_0=\F_{q^{2t_0}}.
	\]
	Let
	\[
	\mathcal C_0=\mathcal S_{\mu_0,\zeta_0}
	=
	\left\{
	\mu_0(v)x+\zeta_0(v)x^{q^{t_0}}
	:
	v\in\LL_0
	\right\}
	\subseteq \End_{\K_0}(\LL_0)
	\]
	be the spread set of a semifield which is two-dimensional over its right
	nucleus $\K_0$. 
	
	Let $h\geq 1$. For each $i=1,\ldots,h$, choose an odd integer $s_i\geq 3$ and
	define recursively
	\[
	t_i=t_{i-1}s_i.
	\]
	Put
	\[
	\K_i=\F_{q^{t_i}},
	\qquad
	\LL_i=\F_{q^{2t_i}}.
	\]
	Thus
	\[
	\K_0\subseteq \K_1\subseteq \cdots \subseteq \K_h,
	\qquad
	\LL_0\subseteq \LL_1\subseteq \cdots \subseteq \LL_h,
	\]
	and, for every $i$, the extension $\LL_i/\K_i$ is quadratic.
	
	For each $i=1,\ldots,h$, let
	\[
	G_i(Y)\in \K_{i-1}[Y]
	\]
	be a monic irreducible polynomial of degree $s_i$, with $G_i(Y)\neq Y$. Let
	$u_i\in\K_i$ be a root of $G_i(Y)$, so that
	\[
	\K_i=\K_{i-1}(u_i).
	\]
	Since $s_i$ is odd, the fields $\LL_{i-1}$ and $\K_i$ are linearly disjoint
	over $\K_{i-1}$. Hence
	\[
	\LL_i=\LL_{i-1}\K_i=\LL_{i-1}(u_i).
	\]
	In particular,
	\[
	\{1,u_i,\ldots,u_i^{s_i-1}\}
	\]
	is an $\LL_{i-1}$-basis of $\LL_i$.
	
	Let
	\[
	\gamma_i:\LL_i\longrightarrow \LL_i,
	\qquad
	x\longmapsto x^{q^{t_i}},
	\]
	be the non-trivial element of $\Gal(\LL_i/\K_i)$. Since $s_i$ is odd, we have
	\[
	t_i=t_{i-1}s_i\equiv t_{i-1}\pmod {2t_{i-1}},
	\]
	and therefore
	\[
	\gamma_i{}_{|\LL_{i-1}}=\gamma_{i-1},
	\]
	where $\gamma_{i-1}$ is the non-trivial element of
	$\Gal(\LL_{i-1}/\K_{i-1})$.
	
	Choose
	\[
	w_i\in\LL_i^\times
	\]
	such that
	\[
	\N_{\LL_i/\K_i}(w_i)=w_i\gamma_i(w_i)=u_i.
	\]
	This is possible because the norm map $\N_{\LL_i/\K_i}:\LL_i^\times\to\K_i^\times$
	is surjective. Also choose
	\[
	\eta_i\in\LL_{i-1}^\times
	\]
	such that
	\[
	\N_{\LL_{i-1}/\K_{i-1}}(\eta_i)=(-1)^{s_i}G_i(0).
	\]
	Again, such an element exists by the surjectivity of the norm map.
	
	We now define recursively additive maps
	\[
	\mu_i,\zeta_i:\LL_i\longrightarrow \LL_i.
	\]
	Assume that $\mu_{i-1},\zeta_{i-1}:\LL_{i-1}\to\LL_{i-1}$ have already been
	defined. Let
	\[
	r_i=\frac{s_i-1}{2}.
	\]
	Every element $y\in\LL_i$ has a unique expression
	\[
	y=v_0+v_1u_i+\cdots+v_{s_i-1}u_i^{s_i-1},
	\qquad
	v_0,\ldots,v_{s_i-1}\in\LL_{i-1}.
	\]
	We define
	\[
	\mu_i(y)
	=
	\mu_{i-1}(v_0)
	+
	\sum_{j=1}^{r_i}v_{2j}u_i^j
	\]
	and
	\[
	\zeta_i(y)
	=
	\left(
	\sum_{j=0}^{r_i-1}v_{2j+1}u_i^j
	+
	\frac{\zeta_{i-1}(v_0)}{\eta_i}u_i^{r_i}
	\right)w_i.
	\]
	The maps $\mu_i$ and $\zeta_i$ are additive, since the expansion of $y$ in the
	$\LL_{i-1}$-basis
	\[
	\{1,u_i,\ldots,u_i^{s_i-1}\}
	\]
	is additive and $\mu_{i-1},\zeta_{i-1}$ are additive.
	
	Define
	\[
	\mathcal C_i
	=
	\mathcal S_{\mu_i,\zeta_i}
	=
	\left\{
	\mu_i(y)x+\zeta_i(y)x^{\gamma_i}
	:
	y\in\LL_i
	\right\}
	\subseteq \End_{\K_i}(\LL_i).
	\]

	\begin{theorem}
		\label{th:iterated-lifting}
		For every $i=0,\ldots,h$, the set $\mathcal C_i$ is the spread set of a
		semifield $\SSS_i$ of order
		\[
		q^{2t_i}=q^{2t s_1\cdots s_i}.
		\]
		Moreover,
		\[
		\NN_r(\SSS_i)\cong \K_i=\F_{q^{t_i}}.
		\]
		In particular, $\SSS_i$ is two-dimensional over its right nucleus.
	\end{theorem}
	
	\begin{proof}
		We argue by induction on $i$. For $i=0$, this is true by assumption.
		
		Assume that $\mathcal C_{i-1}=\mathcal S_{\mu_{i-1},\zeta_{i-1}}$ is the
		spread set of a semifield $\SSS_{i-1}$ of order $q^{2t_{i-1}}$ which is
		two-dimensional over its right nucleus $\K_{i-1}$.
		
		Let
		\[
		R_i=\LL_{i-1}[X_i;\gamma_{i-1}],
		\]
		where $\gamma_{i-1}$ is the non-trivial element of
		$\Gal(\LL_{i-1}/\K_{i-1})$. Since $\gamma_{i-1}$ has order $2$, the center of
		$R_i$ is
		\[
		Z(R_i)=\K_{i-1}[X_i^2].
		\]
		Consider the quotient
		\[
		(R_i)_{G_i}
		=
		\frac{R_i}{R_iG_i(X_i^2)}.
		\]
		The lifting construction of Theorem~\ref{th:liftingspread} gives the spread set
		\[
		\Phi_{G_i,\eta_i}(\mathcal C_{i-1})
		=
		\left\{
		\mu_{i-1}(v_0)
		+
		\sum_{\ell=1}^{s_i-1}v_\ell X_i^\ell
		+
		\frac{\zeta_{i-1}(v_0)}{\eta_i}X_i^{s_i}
		+
		R_iG_i(X_i^2)
		:
		v_0,\ldots,v_{s_i-1}\in\LL_{i-1}
		\right\}
		\subseteq (R_i)_{G_i}.
		\]
		Hence $\Phi_{G_i,\eta_i}(\mathcal C_{i-1})$ is the spread set of a semifield of
		order
		\[
		q^{2t_{i-1}s_i}=q^{2t_i}.
		\]
		By Theorem~\ref{th:right-nucleus-lift}, this semifield has right nucleus
		isomorphic to
		\[
		\K_{i-1}(u_i)=\K_i.
		\]
		
		It remains only to identify this spread set with the explicit semilinear model
		$\mathcal C_i$. We now use the argument of Section \ref{subsec:irreducible-semilinear-maps}. The map
		\[
		\Psi_{w_i}:(R_i)_{G_i}\longrightarrow \End_{\K_i}(\LL_i)
		\]
		is an isomorphism of $\K_i$-algebras. Moreover, if
		$
		T_{w_i}(x)=w_ix^{\gamma_i},
		$
		then
		\[
		T_{w_i}^{2j}(x)=u_i^jx,
		\qquad
		T_{w_i}^{2j+1}(x)=u_i^jw_ix^{\gamma_i}.
		\]
		Therefore, for
		\[
		\overline a
		=
		\mu_{i-1}(v_0)
		+
		\sum_{\ell=1}^{s_i-1}v_\ell X_i^\ell
		+
		\frac{\zeta_{i-1}(v_0)}{\eta_i}X_i^{s_i}
		+
		R_iG_i(X_i^2),
		\]
		we get
		\[
		\begin{array}{rl}
			\Psi_{w_i}(\overline a)(x)
			& =
			\mu_{i-1}(v_0)x
			+
			\sum\limits_{\ell=1}^{s_i-1}v_\ell T_{w_i}^{\ell}(x)
			+
			\frac{\zeta_{i-1}(v_0)}{\eta_i}T_{w_i}^{s_i}(x) \\[6pt]
			& =
			\left(
			\mu_{i-1}(v_0)
			+
			\sum\limits_{j=1}^{r_i}v_{2j}u_i^j
			\right)x \\[6pt]
			& \quad +
			\left(
			\sum\limits_{j=0}^{r_i-1}v_{2j+1}u_i^j
			+
			\frac{\zeta_{i-1}(v_0)}{\eta_i}u_i^{r_i}
			\right)w_ix^{\gamma_i}.
		\end{array}
		\]
		
		Now let
		\[
		y=v_0+v_1u_i+\cdots+v_{s_i-1}u_i^{s_i-1}\in\LL_i.
		\]
		By the definitions of $\mu_i$ and $\zeta_i$, the last expression is precisely
		\[
		\mu_i(y)x+\zeta_i(y)x^{\gamma_i}.
		\]
		Hence
		\[
		\Psi_{w_i}\bigl(\Phi_{G_i,\eta_i}(\mathcal C_{i-1})\bigr)
		=
		\mathcal C_i.
		\]
		Since $\Psi_{w_i}$ is an algebra isomorphism, it preserves nonsingularity and
		the spread-set property. Thus $\mathcal C_i$ is the spread set of a semifield
		$\SSS_i$ of order $q^{2t_i}$, and its right nucleus is isomorphic to $\K_i$.
		This completes the induction.
	\end{proof}
	
	The theorem gives an explicit iterative construction. Starting from a semifield
	$\SSS_0$ of order $q^{2t}$ which is two-dimensional over its right nucleus,
	and choosing odd integers $s_1,\ldots,s_h$ together with irreducible
	polynomials
	\[
	G_i(Y)\in\F_{q^{t_{i-1}}}[Y],
	\qquad
	\deg(G_i)=s_i,
	\]
	one obtains a semifield $\SSS_h$ of order
	\[
	q^{2t s_1\cdots s_h}
	\]
	whose right nucleus is
	\[
	\F_{q^{t s_1\cdots s_h}}.
	\]
	In particular, every step preserves the property of being two-dimensional over
	the right nucleus.

	\section{Explicit constructions}
	\label{sec:explicit-constructions}
	
	In this section we apply the iterative lifting construction developed above to
	the semifields constructed by Ebert, Marino, Polverino and Trombetti
	\cite{ebert2009infinite}.  We use their spread sets as initial spread sets, but
	we adopt the convention of the present paper: spread sets are represented over
	the right nucleus.  Equivalently, we work with the opposite semifields of the
	corresponding EMPT semifields. More precisely, let $\SSS=(\SSS,+,\star)$ be a semifield. The \emph{opposite
	semifield} of $\SSS$ is
	$
	\SSS^{\mathrm{op}}=(\SSS,+,\star^{\mathrm{op}}),
	$
	where
	\[
	a\star^{\mathrm{op}}b:=b\star a
	\qquad (a,b\in \SSS).
	\]
	It has the same multiplicative identity as $\SSS$, and
	its left and right multiplication maps satisfy
	\[
	L_a^{\mathrm{op}}=R_a,
	\qquad
	R_a^{\mathrm{op}}=L_a
	\qquad (a\in \SSS).
	\]
	It follows directly from the definitions that
	\[
	\NN_l(\SSS^{\mathrm{op}})=\NN_r(\SSS),
	\qquad
	\NN_m(\SSS^{\mathrm{op}})=\NN_m(\SSS),
	\qquad
	\NN_r(\SSS^{\mathrm{op}})=\NN_l(\SSS),
	\]
	and
	\[
	Z(\SSS^{\mathrm{op}})=Z(\SSS).
	\] 
	Hence the distinguished nucleus in what follows
	is the right nucleus.
	
	Throughout this section, let $t>1$ be odd and put
	\[
	\K_0=\F_{q^t},
	\qquad
	\LL_0=\F_{q^{2t}},
	\qquad
	\gamma_0:x\longmapsto x^{q^t}.
	\]
	Let $u\in\K_0$ be such that
	\[
	\{1,u,u^2,\ldots,u^{t-1}\}
	\]
	is an $\F_q$-basis of $\K_0$.  Define
	\[
	\mathcal A_u
	=
	\left\langle
	1,u,\ldots,u^{\frac{t-1}{2}}
	\right\rangle_{\F_q},
	\qquad
	\mathcal B_u
	=
	\left\langle
	1,u,\ldots,u^{\frac{t-3}{2}}
	\right\rangle_{\F_q}.
	\]
	Let $a,b,c\in\LL_0^\times$, with $a,c\notin\K_0$, and consider the set
	\[
	S(u,a,b,c)
	=
	\left\{
	x\longmapsto
	(\alpha+a\beta)x+b(\gamma+c\delta)x^{q^t}
	:
	\alpha,\beta\in\mathcal A_u,\ 
	\gamma,\delta\in\mathcal B_u
	\right\}.
	\]
	We assume that the parameters $u,a,b,c$ satisfy one of the admissible
	conditions defining the EMPT families $\mathcal{F}_I,\ldots,\mathcal{F}_V$ in
	\cite{ebert2009infinite}.  Under these assumptions, $S(u,a,b,c)$ is an
	additive spread set of $\K_0$-linear maps of $\LL_0$.
	
	Let
	\[
	\mathcal C_0:=S(u,a,b,c).
	\]
	In the right-nucleus convention used in this paper, $\mathcal C_0$ defines a
	semifield $\SSS_0$ of order
	$
	q^{2t}$.
	Moreover,
	\[
	\NN_r(\SSS_0)\cong \K_0=\F_{q^t},
	\qquad
	\NN_l(\SSS_0)\cong \NN_m(\SSS_0)\cong \F_q,
	\qquad
	Z(\SSS_0)\cong \F_q.
	\]
	Indeed, in the original EMPT convention the semifields in the families
	$\mathcal F_I,\ldots,\mathcal F_V$ are two-dimensional over their left nucleus and have middle
	and right nuclei isomorphic to $\F_q$.  Passing to the opposite semifield
	interchanges the left and right nuclei and leaves the middle nucleus unchanged.

	\subsection{One lifting step}
	
	Let $s\geq 3$ be odd, and write
	\[
	s=2r+1.
	\]
	Let
	\[
	G(Y)\in\K_0[Y]
	\]
	be a monic irreducible polynomial of degree $s$.  Let $\theta\in\F_{q^{ts}}$
	be a root of $G(Y)$, and put
	\[
	\K_1=\K_0(\theta)=\F_{q^{ts}},
	\qquad
	\LL_1=\LL_0(\theta)=\F_{q^{2ts}}.
	\]
	Since $s$ is odd, the fields $\LL_0$ and $\K_1$ are linearly disjoint over
	$\K_0$.  Hence
	\[
	\LL_1=\LL_0\K_1=\LL_0(\theta),
	\]
	and
	\[
	\{1,\theta,\ldots,\theta^{s-1}\}
	\]
	is an $\LL_0$-basis of $\LL_1$.
	
	Choose
	\[
	\eta\in\LL_0^\times
	\]
	such that
	\[
	\N_{\LL_0/\K_0}(\eta)=(-1)^sG(0).
	\]
	Also choose
	\[
	w\in\LL_1^\times
	\]
	such that
	\[
	\N_{\LL_1/\K_1}(w)=ww^{q^{ts}}=\theta.
	\]
	Both choices are possible because norm maps of finite fields are surjective.
	
	Define
	\[
	\mathcal C_1
	\subseteq
	\End_{\K_1}(\LL_1)
	\]
	as the set of all maps
	\[
	x\longmapsto Ax+Bwx^{q^{ts}},
	\]
	where
	\[
	A=
	\alpha+a\beta
	+
	\sum_{j=1}^{r}v_{2j}\theta^j
	\]
	and
	\[
	B=
	\sum_{j=0}^{r-1}v_{2j+1}\theta^j
	+
	\frac{b(\gamma+c\delta)}{\eta}\theta^r,
	\]
	with
	\[
	\alpha,\beta\in\mathcal A_u,
	\qquad
	\gamma,\delta\in\mathcal B_u,
	\qquad
	v_1,\ldots,v_{s-1}\in\LL_0.
	\]
	
	\begin{theorem}
		\label{th:first-explicit-empt-lift}
		The set $\mathcal C_1$ is the spread set of a semifield $\SSS_1$ of order
		$
		q^{2ts}$.
		Moreover,
		\[
		\NN_r(\SSS_1)\cong \K_1=\F_{q^{ts}},
		\]
		and
		\[
		\NN_l(\SSS_1)\cong\F_q,
		\qquad
		\NN_m(\SSS_1)\cong\F_q.
		\]
	\end{theorem}
	
	\begin{proof} The fact that $\mathcal C_1$ is a semifield spread set of order $q^{2ts}$ follows directly from \Cref{th:iterated-lifting}. The same result also gives \[ \NN_r(\SSS_1)\cong \K_1=\F_{q^{ts}}. \] It remains to determine the left and middle nuclei. By \Cref{th:left-middle-nuclei-lift}, the left nucleus of $\SSS_1$ is contained in the left nucleus of the initial semifield $\SSS_0$. Since $\SSS_0$ is the opposite of an EMPT semifield
belonging to one of the families
$\mathcal F_I,\ldots,\mathcal F_V$, its left nucleus
is isomorphic to $\F_q$; see
\cite[Theorem~1.1]{ebert2009infinite}.
Thus, we obtain \[ \F_q\subseteq \NN_l(\SSS_1)\subseteq \NN_l(\SSS_0)\cong \F_q. \] Therefore \[ \NN_l(\SSS_1)\cong \F_q. \] The argument for the middle nucleus is analogous.  \end{proof}
	
	\subsection{The iterated construction}
	
	We now iterate the previous construction by using the argument developed in
	Section~\ref{subsec:iterating-lifting-construction}.
	
	Let $h\geq 1$. For every $i=1,\ldots,h$, let $s_i\geq 3$ be odd and put
	\[
	t_i=ts_1\cdots s_i,
	\qquad
	t_0=t.
	\]
	Set
	\[
	\K_i=\F_{q^{t_i}},
	\qquad
	\LL_i=\F_{q^{2t_i}}.
	\]
	For each $i=1,\ldots,h$, let
	\[
	G_i(Y)\in \K_{i-1}[Y]
	\]
	be a monic irreducible polynomial of degree $s_i$, and let
	$\theta_i\in\K_i$ be a root of $G_i(Y)$. Hence
	\[
	\K_i=\K_{i-1}(\theta_i).
	\]
	Since $s_i$ is odd, the fields $\LL_{i-1}$ and $\K_i$ are linearly disjoint
	over $\K_{i-1}$. Therefore
	\[
	\LL_i=\LL_{i-1}\K_i=\LL_{i-1}(\theta_i),
	\]
	and
	\[
	\{1,\theta_i,\ldots,\theta_i^{s_i-1}\}
	\]
	is an $\LL_{i-1}$-basis of $\LL_i$.
	
	Choose elements
	\[
	\eta_i\in\LL_{i-1}^{\times},
	\qquad
	w_i\in\LL_i^{\times},
	\]
	such that
	\[
	\N_{\LL_{i-1}/\K_{i-1}}(\eta_i)=(-1)^{s_i}G_i(0)
	\]
	and
	\[
	\N_{\LL_i/\K_i}(w_i)=w_iw_i^{q^{t_i}}=\theta_i.
	\]
	These choices are possible by the surjectivity of the norm maps of finite
	fields.
	
	We now describe explicitly the spread set obtained after $h$ lifting steps.
	Let
	\[
	\mathscr P_0=
	\mathcal A_u\times\mathcal A_u\times\mathcal B_u\times\mathcal B_u.
	\]
	For
	\[
	p_0=(\alpha,\beta,\gamma,\delta)\in\mathscr P_0,
	\]
	define
	\[
	A_0(p_0)=\alpha+a\beta,
	\qquad
	B_0(p_0)=b(\gamma+c\delta).
	\]
	Thus
	\[
	\mathcal C_0
	=
	\left\{
	x\longmapsto A_0(p_0)x+B_0(p_0)x^{q^t}
	:
	p_0\in\mathscr P_0
	\right\}.
	\]
	
	We now define the coefficients recursively. For $i=1,\ldots,h$, write
	\[
	r_i=\frac{s_i-1}{2}
	\]
	and put
	\[
	\mathscr P_i
	=
	\mathscr P_{i-1}\times \LL_{i-1}^{s_i-1}.
	\]
	An element of $\mathscr P_i$ will be written as
	\[
	p_i=(p_{i-1};v_{i,1},\ldots,v_{i,s_i-1}),
	\]
	where
	\[
	p_{i-1}\in\mathscr P_{i-1},
	\qquad
	v_{i,1},\ldots,v_{i,s_i-1}\in\LL_{i-1}.
	\]
	Define
	\[
	A_i(p_i)
	=
	A_{i-1}(p_{i-1})
	+
	\sum_{j=1}^{r_i}v_{i,2j}\theta_i^j
	\]
	and
	\[
	B_i(p_i)
	=
	\left(
	\sum_{j=0}^{r_i-1}v_{i,2j+1}\theta_i^j
	+
	\frac{B_{i-1}(p_{i-1})}{\eta_i}\theta_i^{r_i}
	\right)w_i.
	\]
	Then the spread set obtained after $i$ steps is
	\[
	\mathcal C_i
	=
	\left\{
	x\longmapsto A_i(p_i)x+B_i(p_i)x^{q^{t_i}}
	:
	p_i\in\mathscr P_i
	\right\}
	\subseteq \End_{\K_i}(\LL_i).
	\]
	In particular, the final iterated spread set is
	\[
	\mathcal C_h
	=
	\left\{
	x\longmapsto A_h(p_h)x+B_h(p_h)x^{q^{t_h}}
	:
	p_h\in\mathscr P_h
	\right\}
	\subseteq \End_{\K_h}(\LL_h).
	\]

	\begin{theorem}
		\label{th:iterated-explicit-empt-lift}
		The set $\mathcal C_h$ is the spread set of a semifield $\SSS_h$ of order
		\[
		q^{2t_h}=q^{2ts_1\cdots s_h}.
		\]
		Moreover,
		\[
		\NN_r(\SSS_h)\cong \K_h=\F_{q^{ts_1\cdots s_h}},
		\]
		and
		\[
		\NN_l(\SSS_h)\cong \F_q,
		\qquad
		\NN_m(\SSS_h)\cong \F_q, \quad Z(\SSS_h)\cong \F_q.
		\]
	\end{theorem}
	
	\begin{proof}
		The recursive description above is exactly the construction of
		Section~\ref{subsec:iterating-lifting-construction}, applied to the initial
		spread set
		\[
		\mathcal C_0=S(u,a,b,c).
		\]
		Hence, by
		Theorem~\ref{th:iterated-lifting}, $\mathcal C_h$ is a semifield spread set of
		order
		\[
		q^{2t_h}=q^{2ts_1\cdots s_h},
		\]
		and its right nucleus is
		\[
		\NN_r(\SSS_h)\cong \K_h=\F_{q^{t_h}}.
		\]
		
		It remains to determine the left and middle nuclei. At each lifting step,
		Theorem~\ref{th:left-middle-nuclei-lift} shows that the new left and middle
		nuclei are contained in the corresponding nuclei of the previous semifield. Then arguing as in the proof of \Cref{th:first-explicit-empt-lift}, we get the assertion.
		Finally, the center is contained in the intersection of the three nuclei and
		contains the central scalar maps induced by $\F_q$. Thus
		$
		Z(\SSS_h)\cong \F_q$.
	\end{proof}

	\section{The isotopy issue}
	\label{sec:linear-sets-isotopy}
	
	\subsection{Linear sets and semifields} \label{subsec:linearsets}
	In this section we recall the geometric language of linear sets associated
	with spread sets of semifields.  This will be used in the next subsection to
	distinguish the semifields arising from different lifting polynomials.
	
	Let
	\[
	\K=\F_{q^t}
	\qquad\text{and}\qquad
	\LL=\F_{q^{2t}}.
	\]
	We view \(\LL\) as a two-dimensional vector space over \(\K\).  Every
	\(\K\)-linear endomorphism of \(\LL\) has the form
	\[
	\varphi_{A,B}\colon x\longmapsto Ax+Bx^{q^t},
	\qquad
	A,B\in \LL.
	\]
	Thus
	$
	\operatorname{End}_{\K}(\LL)$
	is naturally identified with the \(4\)-dimensional \(\K\)-vector space
	\[
	\LL\times \LL,
	\]
	by sending \(\varphi_{A,B}\) to the pair \((A,B)\).  We shall therefore write
	\[
	\Sigma=\PG(\operatorname{End}_{\K}(\LL),\K)\cong\PG(3,\K).
	\]
	Let \(\mathcal C\subseteq \operatorname{End}_{\K}(\LL)\) be an \(\F_q\)-subspace
	of dimension \(2t\) defining a semifield spread set.  Equivalently,
	\(\mathcal C\) contains the zero map and every nonzero element of
	\(\mathcal C\) is nonsingular. Hence, we may define
	the associated \(\F_q\)-linear set by
	\[
	L(\mathcal C)
	=
	\left\{
	\langle \varphi\rangle_{\K} :
	0\neq \varphi\in\mathcal C
	\right\}
	\subseteq \Sigma.
	\]
	Using the identification \(\operatorname{End}_{\K}(\LL)\simeq \LL\times \LL\), this
	can also be written as
	\[
	L(\mathcal C)
	=
	\left\{
	\langle(A,B)\rangle_{\K} :
	0\neq (A,B)\in\mathcal C
	\right\}
	\subseteq \PG(\LL\times \LL,\K).
	\]
	If \(P \in \Sigma\), the weight of \(P\) with
	respect to \(\mathcal C\) is
	\[
	\operatorname{wt}_{\mathcal C}(P)
	=
	\dim_{\F_q}\bigl(\mathcal C\cap P\bigr).\]
	We shall use repeatedly the following standard criterion, see e.g. \cite[Property 2.3]{polverino2010linear}. Let
	\[
	\ell=\PG(R,\K)
	\]
	be a line of \(\Sigma\), where \(R\) is a two-dimensional \(\K\)-subspace of
	\(\operatorname{End}_{\K}(\LL)\).  Then
	\begin{equation} \label{eq:linecontained}
		\ell\subseteq L(\mathcal C) \Leftrightarrow  \dim_{\F_q}(\mathcal C\cap R)\geq t+1.
	\end{equation}
	
	The nonsingular condition for the maps in \(\mathcal C\) has a useful
	geometric interpretation.  Recall that the map
	$
	\varphi_{A,B}\colon x\longmapsto Ax+Bx^{q^t}$
	is singular if and only if
	$
	\mathrm N_{\LL/\K}(A)=\mathrm N_{\LL/\K}(B)$.
	Consequently, the singular maps define the hyperbolic quadric
	\[
	\mathcal Q
	=
	\left\{
	\langle(A,B)\rangle_{\K} :
	\mathrm N_{\LL/\K}(A)=\mathrm N_{\LL/\K}(B)
	\right\}
	\subseteq \Sigma.
	\]
	Since every nonzero element of a semifield spread set is nonsingular, we have
	\[
	L(\mathcal C)\cap \mathcal Q=\emptyset.
	\]
	The polarity associated with \(\mathcal Q\) will be denoted by
	$
	\perp$.
	In particular, the two coordinate lines
	\[
	\ell=
	\left\{
	\langle(A,0)\rangle_{\K} :
	A\in \LL^\times
	\right\}
	\]
	and
	\[
	\ell^\perp=
	\left\{
	\langle(0,B)\rangle_{\K} :
	B\in \LL^\times
	\right\}
	\]
	are polar to each other.
	
	We shall also use the following standard consequence of the spread-set
	description of isotopy.  If two semifield spread sets
	$
	\mathcal C,\mathcal C'\subseteq\operatorname{End}_{\K}(\LL)$
	define isotopic semifields with the same right nucleus \(\K\), then their
	associated linear sets are equivalent under a semilinear collineation of
	$
	\Sigma=\PG(3,\K)$
	which preserves the hyperbolic quadric \(\mathcal Q\), see \cite[Theorem 2.1]{cardinali2006semifield} and \cite[Theorem 5.3]{ebert2009infinite}.  In particular, such
	a collineation preserves the geometric properties of \(L(\mathcal C)\), such
	as the weights of points, the lines contained in the linear set, and the
	polarity induced by \(\mathcal Q\).
	
	For later use we record the form of the action on the coordinate line \(\ell\), see \cite[Proposition 5.2]{ebert2009infinite}.
	Let \(\Phi\) be a semilinear collineation of \(\Sigma\) preserving
	\(\mathcal Q\), and suppose that
	\[
	\Phi(\ell)=\ell.
	\]
	Then the induced map on \(\ell\) has one of the following two forms:
	\[
	\langle(A,0)\rangle_{\K}
	\longmapsto
	\langle(\lambda A^\rho,0)\rangle_{\K}
	\]
	or
	\[
	\langle(A,0)\rangle_{\K}
	\longmapsto
	\langle(\lambda(A^\rho)^{q^t},0)\rangle_{\K},
	\]
	where
	\[
	\lambda\in \LL^\times
	\qquad\text{and}\qquad
	\rho\in\operatorname{Aut}(\LL).
	\]
	Since $x\mapsto x^{q^t}$ is the nontrivial automorphism
	of $\LL/\K$, both cases can be written, after replacing
	$\rho$ if necessary, in the form
	\[
	\langle(A,0)\rangle_{\K}
	\longmapsto
	\langle(\lambda A^\rho,0)\rangle_{\K}.
	\]
	For further details on linear sets and their connections with semifields, we refer the reader to \cite{cardinali2006semifield,ebert2009infinite,polverino2010linear}.

	We finally recall the field-theoretic lemma that will be used to convert
	equivalence of linear sets into a relation between the roots of the lifting
	polynomials.
	
	\begin{lemma} [see \textnormal{\cite[Lemma~2.4]{Dempwolff2011autotopism}}]
		\label{lem:dempwolff-truncated-subspaces}
		Let \(\K'/\K\) be a field extension of odd degree $s$. Let \(u,w\in \K'\) be such that
		\[
		\K'=\K(u)=\K(w).
		\]
		For $1\leq \ell \leq (s-1)/2$, set
		$
		U=\langle 1,u,u^2,\ldots,u^{\ell}\rangle_{\K}
		$
		and
		$
		W=\langle 1,w,w^2,\ldots,w^{\ell}\rangle_{\K}$.
		If there exists \(\lambda\in \K'^\times\) such that
		$
		U=\lambda W$,
		then
		there exist
		\[
		a,b,c,d\in \K,
		\qquad
		ad-bc\neq0,
		\]
		such that
		\[
		w=\frac{a+bu}{c+du}.
		\]
		Conversely, if \(w\) has this form, then \(U\) is a scalar multiple of
		\(W\).
	\end{lemma}

	\subsection{Isotopy of the lifted semifields}
	\label{sec:isotopy-lifted-semifields}
	
	We now apply the geometric machinery introduced in
	Section~\ref{subsec:linearsets} to the semifields obtained from the
	lifting construction.  We first treat one lifting step and then derive the
	corresponding counting result.
	
	Throughout this section we fix the initial EMPT spread set
	\[
	\mathcal C_0=S(u,a,b,c)
	\subseteq \operatorname{End}_{\K_0}(\LL_0),
	\]
	where
	\[
	\K_0=\F_{q^t},
	\qquad
	\LL_0=\F_{q^{2t}},
	\]
	and \(t\geq 3\) is odd.  We write
	\[
	\mathcal A_u=
	\langle 1,u,\ldots,u^{(t-1)/2}\rangle_{\F_q}
	\]
	and
	\[
	\mathcal B_u=
	\langle 1,u,\ldots,u^{(t-3)/2}\rangle_{\F_q}.
	\]
	Thus
	\[
	\mathcal C_0
	=
	\left\{
	x\longmapsto
	(\alpha+a\beta)x+
	b(\gamma+c\delta)x^{q^t}
	:
	\alpha,\beta\in\mathcal A_u,\ 
	\gamma,\delta\in\mathcal B_u
	\right\}.
	\]
	
	Let
	\[
	s=2r+1\geq 3
	\]
	be odd, and let
	\[
	G(Y)\in \K_0[Y]
	\]
	be monic irreducible of degree \(s\).  Let \(\theta\) be a root of \(G\), and
	put
	\[
	\K_1=\K_0(\theta)=\F_{q^{ts}},
	\qquad
	\LL_1=\LL_0(\theta)=\F_{q^{2ts}}.
	\]
	Choose
	\[
	\eta\in\LL_0^\times
	\qquad\text{and}\qquad
	w\in\LL_1^\times
	\]
	such that
	\[
	\mathrm N_{\LL_0/\K_0}(\eta)=(-1)^sG(0)
	\]
	and
	\[
	\mathrm N_{\LL_1/\K_1}(w)=ww^{q^{ts}}=\theta.
	\]
	
	Let \(\mathcal C_\theta\) be the spread set obtained from
	\(\mathcal C_0\) by the first lifting step.  Hence \(\mathcal C_\theta\)
	consists of the maps
	\[
	x\longmapsto A x+B w x^{q^{ts}},
	\]
	where
	\[
	A=
	\alpha+a\beta+
	\sum_{j=1}^{r}v_{2j}\theta^j
	\]
	and
	\[
	B=
	\sum_{j=0}^{r-1}v_{2j+1}\theta^j
	+
	\frac{b(\gamma+c\delta)}{\eta}\theta^r,
	\]
	with
	\[
	\alpha,\beta\in\mathcal A_u,
	\qquad
	\gamma,\delta\in\mathcal B_u,
	\qquad
	v_1,\ldots,v_{s-1}\in\LL_0.
	\]
	The corresponding semifield has order
	\[
	q^{2ts}
	\]
	and right nucleus isomorphic to
	\[
	\K_1=\F_{q^{ts}}.
	\]
	
	We identify the map
	\[
	x\longmapsto Ax+Bwx^{q^{ts}}
	\]
	with its coefficient pair
	\[
	(A,Bw)\in\LL_1\times\LL_1.
	\]
	Thus the associated linear set is
	\[
	L_\theta=L(\mathcal C_\theta)
	\subseteq\operatorname{PG}(\LL_1\times\LL_1,\K_1)
	\cong\operatorname{PG}(3,\K_1).
	\]
	We denote by
	\[
	\ell=
	\left\{
	\langle(A,0)\rangle_{\K_1}:
	A\in\LL_1^\times
	\right\}
	\]
	the first coordinate line and by
	\[
	\ell^\perp=
	\left\{
	\langle(0,Bw)\rangle_{\K_1}:
	B\in\LL_1^\times
	\right\}
	\]
	the second coordinate line. These two lines are polar
	with respect to the hyperbolic quadric of singular maps.
	
	We first describe the intersections of the defining \(\F_q\)-space of
	\(L_\theta\) with the two coordinate lines.  Put
	\begin{equation} \label{eq:U-theta}
		U_\theta
		=
		\mathcal A_u+a\mathcal A_u
		+
		\sum_{j=1}^{r}\LL_0\theta^j
		\subseteq \LL_1
	\end{equation}
	and
	\begin{equation} \label{eq:yue2}
		W_\theta
		=
		w\left(
		\sum_{j=0}^{r-1}\LL_0\theta^j
		+
		\frac b\eta(\mathcal B_u+c\mathcal B_u)\theta^r
		\right)
		\subseteq \LL_1.
	\end{equation}
	Then
	\[
	\mathcal C_\theta\cap(\LL_1\times\{0\})
	=
	U_\theta\times\{0\}
	\]
	and
	\[
	\mathcal C_\theta\cap(\{0\}\times\LL_1)
	=
	\{0\}\times W_\theta.
	\]
	
	\begin{lemma}
		\label{lem:lift-coordinate-line-weights}
		The line \(\ell\) is contained in \(L_\theta\), whereas \(\ell^\perp\) is not
		contained in \(L_\theta\).  Moreover, every point of \(L_\theta\cap \ell\) has
		weight at most
		\[
		\frac{ts+1}{2},
		\]
		and every point of \(L_\theta\cap \ell^\perp\) has weight at most
		\[
		\frac{ts-1}{2}.
		\]
	\end{lemma}
	
	\begin{proof}
		Since \(a\notin\K_0\), we have
		\[
		\LL_0=\K_0\oplus a\K_0.
		\]
		Therefore
		\[
		U_\theta=P_\theta\oplus aP_\theta,
		\]
		where
		\[
		P_\theta=
		\mathcal A_u+
		\K_0\theta+\K_0\theta^2+\cdots+\K_0\theta^r
		\subseteq\K_1.
		\]
		
		Indeed, for \(j\geq1\),
		\[
		\LL_0\theta^j
		=
		\K_0\theta^j\oplus a\K_0\theta^j.
		\]
		Hence
		\[
		\dim_{\F_q}P_\theta
		=
		\dim_{\F_q}\mathcal A_u+rt
		=
		\frac{t+1}{2}+rt
		=
		\frac{ts+1}{2},
		\]
		and consequently
		\[
		\dim_{\F_q}U_\theta=ts+1.
		\]
		Therefore, by \eqref{eq:linecontained}, we get
		\[
		\ell\subseteq L_\theta.
		\]
		
		We now bound the weights of the points on \(\ell\).  Let
		\[
		P=\langle(\rho,0)\rangle_{\K_1}\in \ell.
		\]
		Write
		\[
		\rho=x+ay,
		\qquad
		x,y\in\K_1,
		\]
		using
		\[
		\LL_1=\K_1\oplus a\K_1.
		\]
		Then
		\[
		\K_1\rho\cap U_\theta
		\]
		is identified with the set of all \(z\in\K_1\) such that
		\[
		z\rho\in P_\theta\oplus aP_\theta.
		\]
		Since
		\[
		z\rho=zx+a(zy),
		\]
		this is equivalent to
		\[
		zx\in P_\theta
		\qquad\text{and}\qquad
		zy\in P_\theta.
		\]
		At least one of \(x,y\) is nonzero.  Therefore the set of such \(z\) has
		\(\F_q\)-dimension at most
		\[
		\dim_{\F_q}P_\theta=\frac{ts+1}{2}.
		\]
		Hence every point of \(L_\theta\cap \ell\) has weight at most
		\[
		\frac{ts+1}{2}.
		\]
		
		We now consider \(\ell^\perp\). 
		Since \(c\notin\K_0\), we have
		\[
		\LL_0=\frac b\eta\K_0\oplus \frac b\eta c\K_0.
		\]
		Thus
		\[
		w^{-1}W_\theta
		=
		\frac b\eta Q_\theta\oplus \frac b\eta cQ_\theta,
		\]
		where
		\[
		Q_\theta=
		\K_0+\K_0\theta+\cdots+\K_0\theta^{r-1}
		+
		\mathcal B_u\theta^r
		\subseteq\K_1.
		\]
		Indeed, for \(0\leq j\leq r-1\),
		\[
		\LL_0\theta^j
		=
		\frac b\eta \K_0\theta^j\oplus \frac b\eta c\K_0\theta^j,
		\]
		and the highest term is
		\[
		\frac b\eta (\mathcal B_u+c\mathcal B_u)\theta^r.
		\]
		Therefore
		\[
		\dim_{\F_q}Q_\theta
		=
		rt+\frac{t-1}{2}
		=
		\frac{ts-1}{2},
		\]
		and so
		\[
		\dim_{\F_q}W_\theta=ts-1.
		\]
		Again \eqref{eq:linecontained} implies that
		\[
		\ell^\perp\not\subseteq L_\theta.
		\]
		
		Arguing as before, we get that every point of \(L_\theta\cap \ell^\perp\) has weight at most
		\[
		\frac{ts-1}{2}. \qedhere
		\]
	\end{proof}
	
	\begin{proposition}
		\label{prop:lift-unique-line}
		The line \(\ell\) is the unique line of
		$
		\PG(\LL_1\times\LL_1,\K_1)$
		contained in \(L_\theta\).
	\end{proposition}
	
	\begin{proof}
		By Lemma~\ref{lem:lift-coordinate-line-weights}, we know that
		$
		\ell \subseteq L_\theta$.
		Suppose, by contradiction, that there exists another line
		\[
		\ell'\neq \ell
		\]
		such that
		$
		\ell'\subseteq L_\theta$.
		
		Let \(R\) and \(S\) be the two-dimensional \(\K_1\)-subspaces of
		\(\LL_1\times\LL_1\) defining \(\ell\) and \(\ell'\), respectively.  Since both
		\(\ell\) and \(\ell'\) are contained in \(L_\theta\), by \eqref{eq:linecontained}, we have
		\[
		\dim_{\F_q}(\mathcal C_\theta\cap R)\geq ts+1
		\]
		and
		\[
		\dim_{\F_q}(\mathcal C_\theta\cap S)\geq ts+1.
		\]
		Since
		$
		\dim_{\F_q}\mathcal C_\theta=2ts,$
		we get
		\[
		(\mathcal C_\theta\cap R)\cap(\mathcal C_\theta\cap S)\neq \{0\}.
		\]
		Therefore \(R\cap S\neq\{0\}\), and the projective lines \(\ell\) and \(\ell'\)
		meet in a point.  Let
		\[
		T=\ell\cap\ell'.
		\]
		
		Let
		\[
		\pi=\langle \ell,\ell'\rangle
		\]
		be the plane generated by \(\ell\) and \(\ell'\), and let \(\Pi\) be the
		corresponding three-dimensional \(\K_1\)-subspace of \(\LL_1\times\LL_1\).
		Then
		\[
		\dim_{\F_q}(\mathcal C_\theta\cap \Pi)
		\geq
		\dim_{\F_q}(\mathcal C_\theta\cap R)
		+
		\dim_{\F_q}(\mathcal C_\theta\cap S)
		-
		\operatorname{wt}_{\mathcal C_\theta}(T).
		\]
		By Lemma~\ref{lem:lift-coordinate-line-weights},
		\[
		\operatorname{wt}_{\mathcal C_\theta}(T)\leq \frac{ts+1}{2}.
		\]
		Therefore
		\[
		\dim_{\F_q}(\mathcal C_\theta\cap \Pi)
		\geq
		(ts+1)+(ts+1)-\frac{ts+1}{2}
		=
		\frac{3ts+3}{2}.
		\]
		
		Now \(\ell^\perp\) meets the plane \(\pi\) in a unique point, say \(P\).  Let
		\(R^\perp\) be the two-dimensional \(\K_1\)-subspace defining \(\ell^\perp\).
		Since
		\[
		\mathcal C_\theta
		=
		(\mathcal C_\theta\cap R)
		+
		(\mathcal C_\theta\cap R^\perp)
		\]
		and \(R\subseteq\Pi\), we have
		\[
		\mathcal C_\theta
		=
		\langle
		\mathcal C_\theta\cap\Pi,\,
		\mathcal C_\theta\cap R^\perp
		\rangle_{\F_q}.
		\]
		Hence
		\[
		2ts
		=
		\dim_{\F_q}\mathcal C_\theta
		\geq
		\dim_{\F_q}(\mathcal C_\theta\cap\Pi)
		+
		\dim_{\F_q}(\mathcal C_\theta\cap R^\perp)
		-
		\operatorname{wt}_{\mathcal C_\theta}(P).
		\]
		Using Lemma~\ref{lem:lift-coordinate-line-weights}, we have
		\[
		\dim_{\F_q}(\mathcal C_\theta\cap R^\perp)=ts-1
		\]
		and
		\[
		\operatorname{wt}_{\mathcal C_\theta}(P)\leq \frac{ts-1}{2}.
		\]
		Thus
		\[
		2ts
		\geq
		\frac{3ts+3}{2}
		+
		(ts-1)
		-
		\frac{ts-1}{2}
		=
		2ts+1,
		\]
		a contradiction.  Therefore \(\ell\) is the unique line contained in
		\(L_\theta\).
	\end{proof}
	
	Let
	$
	H(Y)\in\K_0[Y]$
	be another monic irreducible polynomial of degree \(s\), and let
	$
	\theta'$
	be a root of \(H\).  Let \(\mathcal C_{\theta'}\) be the lifted spread set
	obtained from the same initial spread set \(\mathcal C_0\), using \(\theta'\)
	instead of \(\theta\).
	
	\begin{theorem}
		\label{th:root-isotopy-obstruction}
		Assume that the semifields defined by \(\mathcal C_\theta\) and
		\(\mathcal C_{\theta'}\) are isotopic.  Then there exist
		\[
		\kappa\in\K_0^\times
		\]
		and an automorphism \(\rho\in\operatorname{Aut}(\K_1)\) such that
		\[
		\theta'=\kappa\theta^\rho.
		\]
	\end{theorem}
	
	\begin{proof}
		The two semifields have the same order and their right nuclei are both
		isomorphic to
		$
		\K_1=\F_{q^{ts}}$.
		By the geometric description of isotopy recalled in
		Subsection~\ref{subsec:linearsets}, the isotopy induces a semilinear
		collineation
		$
		\Phi$
		such that
		\[
		\Phi(L_\theta)=L_{\theta'}.
		\]
		
		Moreover, \(\Phi\) preserves the hyperbolic quadric of singular maps. By Proposition~\ref{prop:lift-unique-line}, the line \(\ell\) is the unique line
		contained in \(L_\theta\), and the same is true for \(L_{\theta'}\).  Hence
		\[
		\Phi(\ell)=\ell.
		\]
		By \cite[Proposition 5.2]{ebert2009infinite}, the action induced by \(\Phi\) on \(\ell\) is
		of the form
		\[
		\langle(A,0)\rangle_{\K_1}
		\longmapsto
		\langle(\lambda A^\rho,0)\rangle_{\K_1}
		\]
		for some
		\[
		\lambda\in\LL_1^\times
		\qquad\text{and}\qquad
		\rho\in\operatorname{Aut}(\LL_1).
		\]
		Since finite fields have a unique subfield of each order, \(\rho\) stabilizes
		\(\K_1\), \(\LL_0\), and \(\K_0\).  We shall denote again by \(\rho\) its
		restriction to \(\K_1\). 
		
		We have
		\[
		\LL_1=\LL_0(\theta^\rho)=\LL_0(\theta'),
		\]
		because \(\LL_0\) and \(\K_1\) are linearly disjoint over \(\K_0\). The collineation \(\Phi\) maps the defining \(\F_q\)-subspace of
		\(L_\theta\) onto the defining \(\F_q\)-subspace of \(L_{\theta'}\). Since
		\(\Phi(\ell)=\ell\), it maps the intersection with the first coordinate space onto
		the corresponding intersection for \(\theta'\). Therefore
		\[
		\lambda U_\theta^\rho=U_{\theta'}.
		\]
		For \(z\in\{\theta^\rho,\theta'\}\), set
		\[
		S_z=\LL_0+\LL_0z+\cdots+\LL_0z^r,
		\]
		\[
		J_z=\LL_0z+\cdots+\LL_0z^r,
		\]
		and
		\[
		R_z=\LL_0+\LL_0z+\cdots+\LL_0z^{r-1}.
		\]
		
		Recalling the definition of $U_\theta$ in \eqref{eq:U-theta},
		we claim that $J_{\theta'}$ is the largest $\LL_0$-subspace contained
		in $U_{\theta'}$, and similarly $J_{\theta^\rho}$ is the largest
		$\LL_0$-subspace contained in $U_\theta^\rho$. Indeed, the $\LL_0$-component of an element of $U_{\theta'}$,
		with respect to the decomposition $S_{\theta'}=\LL_0\oplus J_{\theta'}$,
		lies in
		\[
		\mathcal A_u+a\mathcal A_u,
		\]
		which has \(\F_q\)-dimension \(t+1\).  Since \(t\geq3\), this is strictly
		smaller than
		\[
		\dim_{\F_q}\LL_0=2t.
		\]
		Hence
		any \(\LL_0\)-subspace contained in \(U_{\theta'}\) must have zero $\LL_0$-component, and is therefore contained in \(J_{\theta'}\).  The same argument applies
		to \(U_\theta^\rho\), because \(\rho\) stabilizes \(\LL_0\) and preserves
		\(\F_q\)-dimensions.  Consequently, from \(\lambda U_\theta^\rho=U_{\theta'}\)
		we get
		\[
		\lambda S_{\theta^\rho}=S_{\theta'}
		\qquad\text{and}\qquad
		\lambda J_{\theta^\rho}=J_{\theta'}.
		\]
		
		Applying Lemma~\ref{lem:dempwolff-truncated-subspaces} to the equality
		\(\lambda S_{\theta^\rho}=S_{\theta'}\), with base field \(\LL_0\), gives
		\[
\theta'=\frac{A+B\theta^\rho}{C+D\theta^\rho}
		\]
		for some
		\[
		A,B,C,D\in\LL_0,
		\qquad
		AD-BC\neq0.
		\]
		We now show that \(A=0\).  Put
		\[
		\Delta=(C+D\theta^\rho)^r.
		\]
		Then
		\begin{equation} \label{eq:yue3}
			\Delta S_{\theta'}=S_{\theta^\rho}.
		\end{equation}
		Since \(\lambda S_{\theta^\rho}=S_{\theta'}\), we have
		\[
		\Delta\lambda S_{\theta^\rho}=S_{\theta^\rho}.
		\]
		We shall use the following elementary observation: if \(\alpha S_{\theta^\rho}=S_{\theta^\rho}\) for
		some \(\alpha\in\LL_1^\times\), then \(\alpha\in\LL_0^\times\).  Indeed,
		\(S_{\theta^\rho}\) is then a vector space over \(\LL_0(\alpha)\). Hence
		\([\LL_0(\alpha):\LL_0]\) divides \(\dim_{\LL_0}S_{\theta^\rho}=r+1\).  On the other hand,
		\([\LL_0(\alpha):\LL_0]\) divides \([\LL_1:\LL_0]=2r+1\). Since
		\[
		\gcd(r+1,2r+1)=1,
		\]
		we get \(\alpha\in\LL_0^\times\). It follows that
		\begin{equation} \label{eq:yue4}
			\Delta\lambda\in\LL_0^\times.
		\end{equation}
		Therefore
		\[
		\Delta J_{\theta'}
		=
		\Delta\lambda J_{\theta^\rho}
		=
		J_{\theta^\rho}.
		\]
		In particular,
		\[
		\Delta (\theta')^r\in J_{\theta^\rho}.
		\]
		But
		\[
		\Delta (\theta')^r
		=
		(C+D\theta^\rho)^r\left(\frac{A+B\theta^\rho}{C+D\theta^\rho}\right)^r
		=
		(A+B\theta^\rho)^r.
		\]
		Since
		\[
		J_{\theta^\rho}=\LL_0\theta^\rho+\cdots+\LL_0(\theta^\rho)^r
		\]
		consists of the elements of \(S_{\theta^{\rho}}\) with zero constant coefficient, the constant
		coefficient of \((A+B\theta^\rho)^r\) must be zero.  This coefficient is \(A^r\), and so
		$
		A=0$.
		Thus
		\[
		\theta'=\frac{B\theta^\rho}{C+D\theta^\rho},
		\qquad
		BC\neq0.
		\]
		
		It remains to prove that \(D=0\). Since \(\Phi\) preserves the hyperbolic
		quadric and \(\Phi(\ell)=\ell\), it also preserves the polar line \(\ell^\perp\). Hence
		\[
		\Phi(\ell^\perp)=\ell^\perp.
		\]
		Using again the form of the action on a coordinate line, there exists
		\(\mu\in\LL_1^\times\) such that the induced action on \(\ell^\perp\) gives
		\[
		\mu W_\theta^\rho=W_{\theta'}.
		\]
		Write
		\[
		W_\theta=wT_\theta,
		\]
		where
		\[
		T_\theta=
		\sum_{j=0}^{r-1}\LL_0\theta^j
		+
		\frac b\eta(\mathcal B_u+c\mathcal B_u)\theta^r.
		\]
		The \(\LL_0\)-span of \(T_{\theta'}\) is \(S_{\theta'}\), and the largest
		\(\LL_0\)-subspace contained in \(T_{\theta'}\) is \(R_{\theta'}\). Indeed, the
		coefficient of \({\theta'}^r\) belongs to
		\[
		\frac b{\eta'}(\mathcal B_u+c\mathcal B_u),
		\]
		which has \(\F_q\)-dimension \(t-1<2t\).  The same argument applies to \(T_\theta^\rho\).  Therefore,
		after removing the nonzero scalar factors \(w^\rho\) and the corresponding
		choice \(w'\) for \(\theta'\), the equality \(\mu W_\theta^\rho=W_{\theta'}\)
		gives a scalar \(\nu\in\LL_1^\times\) such that
		\[
		\nu S_{\theta^\rho}=S_{\theta'}
		\qquad\text{and}\qquad
		\nu R_{\theta^\rho}=R_{\theta'}.
		\]
		By \eqref{eq:yue3}, we obtain
		\[
		\Delta\nu S_{\theta^\rho}=S_{\theta^\rho}.
		\]
		By the same computation for \eqref{eq:yue4},
		\[
		\Delta\nu\in\LL_0^\times.
		\]
		
		Therefore
		\[
		\Delta R_{\theta'}
		=
		\Delta\nu R_{\theta^\rho}
		=
		R_{\theta^\rho}.
		\]
		Since \(1\in R_{\theta'}\), this implies
		\[
		\Delta=(C+D\theta^\rho)^r\in R_{\theta^\rho}.
		\]
		But
		\[
		R_{\theta^\rho}=\LL_0+\LL_0\theta^\rho+\cdots+\LL_0(\theta^\rho)^{r-1}
		\]
		consists of the elements of \(S_{\theta^\rho}\) with zero coefficient of \((\theta^\rho)^r\). The
		coefficient of \((\theta^\rho)^r\) in \((C+D\theta^\rho)^r\) is \(D^r\). Hence
		$
		D=0$.
		Thus
		\[
		\theta'=\frac BC \theta^\rho.
		\]
		Set
		\[
		\kappa=\frac BC\in\LL_0^\times.
		\]
		Since \(\theta^\rho,\theta'\in\K_1\), we also have
		\[
		\kappa=\frac{\theta'}{\theta^\rho}\in\K_1^\times.
		\]
		Finally,
		\[
		\K_1\cap\LL_0=\K_0,
		\]
		because \([\K_1:\K_0]=s\) is odd and \([\LL_0:\K_0]=2\). Hence
		\[
		\kappa\in\K_0^\times.
		\]
		This proves the theorem.
	\end{proof}

	We now count the number of isotopy classes obtained by varying the lifting
	polynomial \(G\). Let
	\[
	Q=|\K_0|=q^t.
	\] Denote by
	\[
	I_Q(s)=\frac{1}{s}\sum_{d\mid s}\mu(d)Q^{s/d}
	\]
	the number of monic irreducible polynomials of degree
	$
	s$
	over
	$
	\K_0$.

	\begin{theorem}
		\label{th:first-lift-count}
		Fix the initial spread set
		$
		\mathcal C_0=S(u,a,b,c)$
		and fix an odd integer
		$
		s\geq3$. Let
		$
		\mathfrak F_s$
		be the family of semifields arising from the first lifting step. Then
		$
		\mathfrak F_s$
		contains at least
		\[
		\left\lceil
		\frac{I_Q(s)}
		{|\operatorname{Aut}(\K_0)|(Q-1)}
		\right\rceil
		\]
		pairwise non-isotopic semifields of order
		\[
		q^{2ts}=Q^{2s}.
		\]
	\end{theorem}
	
	\begin{proof}
		
		By Theorem~\ref{th:root-isotopy-obstruction}, two such polynomials can give
		isotopic semifields only if their roots are related by
		\[
		\theta'=\kappa\theta^\rho,
		\]
		with
		\[
		\kappa\in\K_0^\times
		\qquad\text{and}\qquad
		\rho\in\operatorname{Aut}(\K_1).
		\]
		
		Let \(\tau=\rho|_{\K_0}\in\operatorname{Aut}(\K_0)\). If \(G(Y)\) is the
		minimal polynomial of \(\theta\) over \(\K_0\), then \(\theta^\rho\) is a root
		of the polynomial obtained by applying \(\tau\) to the coefficients of \(G\),
		which we denote by \(G^\tau(Y)\). Therefore \(\kappa\theta^\rho\) is a root of
		\[
		G_{\kappa,\tau}(Y)
		=
		\kappa^sG^\tau\left(\frac Y\kappa\right).
		\]
		For a fixed \(G\), the number of polynomials that can arise in this way is at
		most
		\[
		|\K_0^\times|\,|\operatorname{Aut}(\K_0)|
		=
		(Q-1)|\operatorname{Aut}(\K_0)|.
		\]
		The automorphisms of \(\K_1\) which fix \(\K_0\) only replace a root by another
		root of the same irreducible polynomial, and hence do not produce further
		polynomials.  Therefore each isotopy class contains at most
		\[
		(Q-1)|\operatorname{Aut}(\K_0)|
		\]
		choices of the polynomial \(G\). It follows that the number of isotopy classes
		is at least
		\[
		\left\lceil
		\frac{I_Q(s)}{|\operatorname{Aut}(\K_0)|(Q-1)}
		\right\rceil,
		\]
		as claimed.
	\end{proof}
	
	Theorem~\ref{th:first-lift-count} shows that, for fixed odd integers
	\(t,s\geq 3\) and \(q=p^f\), where \(p\) is prime, the first lifting step gives at least
	\[
	\Omega_{t,s}\!\left(
	\frac{q^{t(s-1)}}{sf}
	\right)
	\]
	pairwise non-isotopic semifields of order
	\[
	N=q^{2ts},
	\]
	which are two-dimensional over their right nucleus of order \(q^{ts}\), as
	\(q\to+\infty\).

	\subsection{Isotopy for the iteration}
	\label{subsec:isotopy-iteration}
	
	We now explain how the isotopy criterion obtained for the first lifting step
	extends to every level of the iterative construction.  The point is that, at
	each step, the new spread set has the same coordinate-line structure as in the
	first lift.  Hence one can vary the polynomial used in the quotient
	construction at any prescribed level and still obtain many pairwise
	non-isotopic semifields.
	
	We keep the notation of Section~\ref{sec:explicit-constructions}.  Thus
	\[
	\K_i=\F_{q^{t_i}},
	\qquad
	\LL_i=\F_{q^{2t_i}},
	\qquad
	t_i=ts_1\cdots s_i,
	\]
	and
	\[
	\mathcal C_i
	=
	\left\{
	x\longmapsto A_i(p_i)x+B_i(p_i)x^{q^{t_i}}
	:
	p_i\in\mathscr P_i
	\right\}
	\subseteq \End_{\K_i}(\LL_i)
	\]
	is the spread set obtained after \(i\) lifting steps.
	
	Let
	\[
	\ell_i=
	\{\langle(A,0)\rangle_{\K_i}:A\in\LL_i^\times\}
	\]
	and
	\[
	\ell_i^\perp=
	\{\langle(0,B)\rangle_{\K_i}:B\in\LL_i^\times\}
	\]
	be the two coordinate lines in
	\[
	\PG(\LL_i\times\LL_i,\K_i)\cong \PG(3,\K_i).
	\]
	We denote by
	\[
	U_i=
	\{A\in\LL_i:(A,0)\in \mathcal C_i\}
	\]
	and
	\[
	W_i=
	\{B\in\LL_i:(0,B)\in \mathcal C_i\}.
	\]
	
	We first record the recursive form of these two coordinate sections.
	
	\begin{lemma}
		\label{lem:coordinate-sections-iteration}
		For every \(i\geq 0\), there exist \(\F_q\)-subspaces
		\[
		\mathcal P_i,\mathcal Q_i\subseteq \K_i
		\]
		and an element \(d_i\in\LL_i^\times\) such that
		\[
		U_i=\mathcal P_i\oplus a \mathcal P_i
		\]
		and
		\[
		W_i=d_i(\mathcal Q_i\oplus c\mathcal Q_i).
		\]
		Moreover,
		\[
		\dim_{\F_q}\mathcal P_i=\frac{t_i+1}{2},
		\qquad
		\dim_{\F_q}\mathcal Q_i=\frac{t_i-1}{2}.
		\]
	\end{lemma}
	
	\begin{proof}
		We argue by induction on $i$, extending the decomposition
		used in the proof of Lemma~\ref{lem:lift-coordinate-line-weights} to every step of the iteration.
		For $i=0$, we have
		\[
		U_0=\mathcal A_u+a\mathcal A_u
		\]
		and
		\[
		W_0=b(\mathcal B_u+c\mathcal B_u).
		\]
		Hence we may take
		\[
		\mathcal P_0=\mathcal A_u,
		\qquad
		\mathcal Q_0=\mathcal B_u,
		\qquad
		d_0=b.
		\]
		The asserted dimensions are precisely
		\[
		\dim_{\F_q}\mathcal A_u=\frac{t+1}{2},
		\qquad
		\dim_{\F_q}\mathcal B_u=\frac{t-1}{2}.
		\]
		
		Assume now that the assertion holds for \(i-1\).  Write
		\[
		s_i=2r_i+1.
		\]
		By the recursive definition of \(\mathcal C_i\), an element of
		\(\mathcal C_i\) has first coordinate
		\[
		A_i
		=
		A_{i-1}
		+
		\sum_{j=1}^{r_i}v_{i,2j}\theta_i^j
		\]
		and second coordinate
		\[
		B_i
		=
		\left(
		\sum_{j=0}^{r_i-1}v_{i,2j+1}\theta_i^j
		+
		\frac{B_{i-1}}{\eta_i}\theta_i^{r_i}
		\right)w_i.
		\]
		Therefore
		\[
		U_i
		=
		U_{i-1}
		+
		\LL_{i-1}\theta_i
		+\cdots+
		\LL_{i-1}\theta_i^{r_i}.
		\]
		Since \(a\notin\K_0\), and since
		\[
		\K_{i-1}\cap\LL_0=\K_0,
		\]
		we have \(a\notin\K_{i-1}\). Hence
		\[
		\LL_{i-1}=\K_{i-1}\oplus a\K_{i-1}.
		\]
		Using the induction hypothesis
		\[
		U_{i-1}=\mathcal P_{i-1}\oplus a\mathcal P_{i-1},
		\]
		we get
		\[
		U_i=\mathcal P_i\oplus a\mathcal P_i,
		\]
		where
		\[
		\mathcal P_i=
		\mathcal P_{i-1}
		+
		\K_{i-1}\theta_i
		+\cdots+
		\K_{i-1}\theta_i^{r_i}
		\subseteq \K_i.
		\]
		Thus
		\[
		\dim_{\F_q}\mathcal P_i
		=
		\dim_{\F_q}\mathcal P_{i-1}+r_it_{i-1}
		=
		\frac{t_{i-1}+1}{2}
		+
		\frac{s_i-1}{2}t_{i-1}
		=
		\frac{t_i+1}{2}.
		\]
		
		Similarly,
		\[
		W_i
		=
		w_i\left(
		\LL_{i-1}
		+
		\LL_{i-1}\theta_i
		+\cdots+
		\LL_{i-1}\theta_i^{r_i-1}
		+
		\frac{W_{i-1}}{\eta_i}\theta_i^{r_i}
		\right).
		\]
		By the induction hypothesis,
		\[
		W_{i-1}=d_{i-1}(\mathcal Q_{i-1}\oplus c\mathcal Q_{i-1}).
		\]
		Set
		\[
		d_i=w_i\eta_i^{-1}d_{i-1}.
		\]
		Since \(c\notin\K_{i-1}\), we have
		\[
		\LL_{i-1}=\K_{i-1}\oplus c\K_{i-1}.
		\]
		Therefore
		\[
		W_i=d_i(\mathcal Q_i\oplus c\mathcal Q_i),
		\]
		where
		\[
		\mathcal Q_i=
		\K_{i-1}
		+
		\K_{i-1}\theta_i
		+\cdots+
		\K_{i-1}\theta_i^{r_i-1}
		+
		\mathcal Q_{i-1}\theta_i^{r_i}
		\subseteq \K_i.
		\]
		Hence
		\[
		\dim_{\F_q}\mathcal Q_i
		=
		r_it_{i-1}
		+
		\dim_{\F_q}\mathcal Q_{i-1}
		=
		\frac{s_i-1}{2}t_{i-1}
		+
		\frac{t_{i-1}-1}{2}
		=
		\frac{t_i-1}{2}.
		\]
		This completes the induction.
	\end{proof}
	
	\begin{corollary}
		\label{cor:unique-line-iteration}
		For every \(i\geq 1\), the line \(\ell_i\) is the unique line of
		\[
		\PG(\LL_i\times\LL_i,\K_i)
		\]
		contained in the linear set \(L(\mathcal C_i)\).
	\end{corollary}
	
	\begin{proof}
		By Lemma~\ref{lem:coordinate-sections-iteration}, we have
		\[
		\dim_{\F_q}U_i=t_i+1.
		\]
		Hence, by \eqref{eq:linecontained}, we have
		\[
		\ell_i\subseteq L(\mathcal C_i).
		\]
		Moreover, every point of \(L(\mathcal C_i)\cap \ell_i\) has weight at most
		\[
		\dim_{\F_q}\mathcal P_i=\frac{t_i+1}{2}.
		\]
		Indeed, write
		\[
		\rho=x+ay,\qquad x,y\in\K_i,
		\]
		where $\rho\ne0$. Since $U_i=\mathcal P_i\oplus a\mathcal P_i$, we have
		\[
		\K_i\rho\cap U_i
		=
		\{z\rho:z\in\K_i,\ zx\in \mathcal P_i,\ zy\in \mathcal P_i\}.
		\]
		If $x\ne0$, the map
		\[
		\K_i\rho\cap U_i\longrightarrow \mathcal P_i,
		\qquad z\rho\longmapsto zx,
		\]
		is $\mathbb F_q$-linear and injective. If $x=0$, then
		$y\ne0$, and the same conclusion holds using
		$z\rho\longmapsto zy$ instead. Therefore
		\[
		\dim_{\mathbb F_q}(\K_i\rho\cap U_i)
		\le \dim_{\mathbb F_q}\mathcal P_i
		=\frac{t_i+1}{2}.
		\]
		The left-hand side is precisely the weight of
		$\langle(\rho,0)\rangle_{\K_i}$ in $L(\C_i)$,
		which proves the claimed bound.
		
		Similarly,
		\[
		\dim_{\F_q}W_i=t_i-1,
		\]
		so \(\ell_i^\perp\) is not contained in \(L(\mathcal C_i)\), and every point of
		\(L(\mathcal C_i)\cap \ell_i^\perp\) has weight at most
		\[
		\dim_{\F_q}\mathcal Q_i=\frac{t_i-1}{2}.
		\]
		The same dimension argument used in Proposition~\ref{prop:lift-unique-line}
		now applies verbatim and shows that no line other than \(\ell_i\) can be contained
		in \(L(\mathcal C_i)\).
	\end{proof}
	
	We can now repeat the root-obstruction argument at an arbitrary step of the
	iteration.
	
	\begin{theorem}
		\label{th:relative-root-obstruction}
		Fix \(i\geq 1\), and fix the spread set \(\mathcal C_{i-1}\) obtained after
		the first \(i-1\) steps.  Let
		\[
		s_i=2r_i+1\geq 3
		\]
		be odd.  Let
		\[
		G_i(Y),H_i(Y)\in\K_{i-1}[Y]
		\]
		be monic irreducible polynomials of degree \(s_i\).  Let \(\theta_i\) be a
		root of \(G_i\), and let \(\theta_i'\) be a root of \(H_i\).  Let
		\[
		\mathcal C_i(\theta_i)
		\qquad\text{and}\qquad
		\mathcal C_i(\theta_i')
		\]
		be the spread sets obtained from \(\mathcal C_{i-1}\) by the \(i\)-th lifting
		step using \(\theta_i\) and \(\theta_i'\), respectively.
		
		If the corresponding semifields are isotopic, then there exist
		\[
		\kappa\in\K_{i-1}^{\times}
		\]
		and
		\[
		\rho\in\operatorname{Aut}(\K_i)
		\]
		such that
		\[
		\theta_i'=\kappa\theta_i^\rho.
		\]
	\end{theorem}
	
	\begin{proof}
		The proof is the same as that of
		Theorem~\ref{th:root-isotopy-obstruction}, with
		\[
		\K_0,\LL_0,\K_1,\LL_1
		\]
		replaced by
		\[
		\K_{i-1},\LL_{i-1},\K_i,\LL_i.
		\]
	\end{proof}
	
	We next show that the liftings of non-isotopic semifields remain non-isotopic.
		
		\begin{theorem}
			\label{thm:equivalence}
			Fix \(i\geq 1\), and put
			$\K=\K_{i-1}$, $\LL=\LL_{i-1}$, $\K'=\K_i$, $\LL'=\LL_i$.
			Let further $s=s_i=2r+1\geq3$, and $\mathcal C,\mathcal C'
			\subseteq
			\End_{\K}(\LL)$
			be two spread sets occurring after $i-1$ steps of the construction, possibly
			arising from different choices at the preceding steps. Choose monic
			irreducible polynomials
			$G(Y),H(Y)\in\K[Y]$
			of degree $s$, roots
			$\theta,\theta'\in\K',$
			and admissible elements
			$\eta,\eta'\in\LL^\times,
			w,w'\in{\LL'}^\times,$
			so that
			$\N_{\LL/\K}(\eta)=(-1)^sG(0),$
			$\N_{\LL/\K}(\eta')=(-1)^sH(0),$ 
			$\N_{\LL'/\K'}(w)=\theta,$
			$\N_{\LL'/\K'}(w')=\theta'.$
			
			Let
			$\widetilde{\mathcal C}_{\theta}$ and
			$\widetilde{\mathcal C'}_{\theta'}$
			be the corresponding lifted spread sets. If the semifields defined by these
			two lifted spread sets are isotopic, then the semifields defined by
			$\mathcal C$ and $\mathcal C'$ are isotopic.
		\end{theorem}
		\begin{proof}
			We identify the map
			\[
			x\mapsto Ax+Bx^{q^{t_{i-1}}}
			\]
			with the pair \((A,B)\in\LL\times\LL\). In this notation the first lifted
			spread set has the form
			\[
			\widetilde{\mathcal C}_{\theta}
			=
			\left\{
			\left(
			A+\sum_{j=1}^{r}a_j\theta^j,
			\left(
			\sum_{j=0}^{r-1}b_j\theta^j
			+\frac{B}{\eta}\theta^r
			\right)w
			\right):
			(A,B)\in\mathcal C,\ 
			a_j,b_j\in\LL
			\right\},
			\]
			and we can find a similar representation for
			\(\widetilde{\mathcal C'}_{\theta'}\).
			
			Suppose that the lifted semifields are isotopic. By
			Corollary~\ref{cor:unique-line-iteration}, the first coordinate line is the
			unique line contained in each of the associated linear sets. Hence the
			semilinear collineation induced by the isotopy fixes this line.
			Using the spread-set form of isotopy from
			Proposition~\ref{prop:isotopy-equivalence-right-nucleus}, and absorbing the
			non-trivial automorphism of \(\LL'/\K'\) into the field automorphism if
			necessary, the action on coefficient pairs is then of the form
			\begin{equation}
				\label{eq:parent-recovery-diagonal-action}
				(A,B)
				\mapsto
				(\lambda A^\rho,\mu B^\rho)
			\end{equation}
			for some
			$\lambda,\mu\in{\LL'}^\times$ and $\rho\in\Aut(\LL')$
			with
			\begin{equation}
				\label{eq:parent-recovery-big-norms}
				\N_{\LL'/\K'}(\lambda)
				=
				\N_{\LL'/\K'}(\mu).
			\end{equation}
			For convenience, put
			\[
			\tau=\rho|_{\K}, \;
			x=\theta^\rho, \;
			y=\theta'.
			\]
			The argument in the proof of
			Theorem~\ref{th:root-isotopy-obstruction} applies also when the two original
			spread sets are different. Indeed, by
			Lemma~\ref{lem:coordinate-sections-iteration}, the first coordinate section of
			each original spread set has $\F_q$-dimension $t_{i-1}+1<2t_{i-1},$ whereas the second
			coordinate section has dimension $t_{i-1}-1<2t_{i-1}$. These are the only
			properties of the spread set used to identify the largest $\LL$-subspaces in the
			two lifted coordinate sections. Applying
			Lemma~\ref{lem:dempwolff-truncated-subspaces} to their $\LL$-spans, and then
			using both coordinate sections as in the proof of
			Theorem~\ref{th:root-isotopy-obstruction}, gives
			\begin{equation}
				\label{eq:parent-recovery-root-relation}
				y=\kappa x
			\end{equation}
			for some
			\[
			\kappa\in\K^\times.
			\]
			
			For \(z\in\{x,y\}\), put
			\[
			S_z=\LL+\LL z+\cdots+\LL z^r.
			\]
			The equality of the \(\LL\)-spans of the first coordinate sections gives
			\[
			\lambda S_x=S_y.
			\]
			Since \(y=\kappa x\) with \(\kappa\in\K^\times\), we have \(S_y=S_x\).
			The stabilizer computation used in the proof of
			Theorem~\ref{th:root-isotopy-obstruction} therefore yields that
			$\lambda\in\LL^\times$.
			
			Similarly, after removing the factors $w^\rho$ and $w'$ from the second
			coordinate sections, put $\nu=\frac{\mu w^\rho}{w'}$.
			Their $\LL$-spans then satisfy
			$\nu S_x=S_y=S_x$,
			and we conclude that $\nu\in\LL^\times$.
			
			We now pass from the lifted spread sets to the original spread sets. Since
			$y=\kappa x$, the free-coefficient subspaces in the two lifts agree after
			the transformation in Eq.~\eqref{eq:parent-recovery-diagonal-action}. Taking the
			quotient by these free coefficients, or equivalently comparing the constant
			coefficient in the first coordinate and the coefficient of $x^r$ in the
			second coordinate, gives
			\begin{equation}
				\label{eq:parent-recovery-parent-equivalence}
				\mathcal C'
				=
				\left\{
				(\lambda A^\rho,\beta B^\rho):
				(A,B)\in\mathcal C
				\right\},
			\end{equation}
			where
			\begin{equation}
				\label{eq:parent-recovery-beta}
				\beta
				=
				\nu\frac{\eta'}{\eta^\rho}\kappa^{-r}
				\in\LL^\times.
			\end{equation}
			
			It remains to show that the two coordinate factors in
			Eq.~\eqref{eq:parent-recovery-parent-equivalence} have the same norm over
			$\K$. 
			By Eqs.~\eqref{eq:parent-recovery-big-norms},~\eqref{eq:parent-recovery-root-relation} and the definitions of $\nu,x,y$, we have
			\[
			\N_{\LL/\K}(\nu)
			=
			\N_{\LL'/\K'}(\nu)
			=
			\N_{\LL'/\K'}(\mu)\frac{x}{y}
			=
			\frac{\N_{\LL/\K}(\lambda)}{\kappa}.
			\]
			Moreover, \(x\) is a root of \(G^\tau\), while \(y=\kappa x\) is a root of
			\(H\). Since both polynomials are monic and irreducible of degree \(s\),
			\[
			H(Y)=\kappa^sG^\tau\left(\frac{Y}{\kappa}\right).
			\]
			Consequently,
			\[
			\frac{\N_{\LL/\K}(\eta')}{\N_{\LL/\K}(\eta^\rho)}
			=
			\frac{H(0)}{G^\tau(0)}
			=
			\kappa^s.
			\]
			Since \(\kappa\in\K\) and \(\LL/\K\) is quadratic, we also have
			\[
			\N_{\LL/\K}(\kappa^{-r})=\kappa^{-2r}.
			\]
			Using Eq.~\eqref{eq:parent-recovery-beta} and \(s=2r+1\), we obtain
			\begin{align*}
				\N_{\LL/\K}(\beta)
				&=
				\N_{\LL/\K}(\nu)
				\frac{\N_{\LL/\K}(\eta')}{\N_{\LL/\K}(\eta^\rho)}
				\N_{\LL/\K}(\kappa^{-r})\\
				&=
				\frac{\N_{\LL/\K}(\lambda)}{\kappa}
				\kappa^s
				\kappa^{-2r}
				=
				\N_{\LL/\K}(\lambda).
			\end{align*}
			
			Thus $\N_{\LL/\K}(\beta/\lambda)=1$. Let $\gamma$ be the non-trivial
			automorphism of $\LL/\K$. So there exists
			$d\in\LL^\times$ such that
			$\beta/\lambda=d^{\gamma-1}$.
			Put $c=\lambda/d$. This means
			$\lambda=cd,
			\beta=cd^\gamma.$
			
			If
			\[
			T_{A,B}:z\mapsto Az+Bz^\gamma,
			\]
			and \(M_c,M_d\) denote scalar multiplication by $c,d$, respectively, then
			\[
			M_c\circ T_{A^\rho,B^\rho}\circ M_d
			=
			T_{\lambda A^\rho,\beta B^\rho}.
			\]
			Now, it follows from Eq.~\eqref{eq:parent-recovery-parent-equivalence} and
			Proposition~\ref{prop:isotopy-equivalence-right-nucleus} that
			$\mathcal C$ and $\mathcal C'$ define isotopic semifields.
		\end{proof}
		
		In particular, 
		Theorem~\ref{thm:equivalence} shows that lifts of non-isotopic 
		semifields are non-isotopic.

	We now give a concrete consequence of the previous criterion.  The idea is
	that, at each level of the iteration, once the previous semifield has been
	fixed, the polynomial used in the quotient construction can be varied in many
	essentially different ways.
	
	Let
	\[
	Q_{i-1}=|\K_{i-1}|=q^{t_{i-1}}.
	\]
	For an odd integer \(s_i\geq 3\), denote by
	\[
	I_{Q_{i-1}}(s_i)
	=
	\frac{1}{s_i}
	\sum_{d\mid s_i}\mu(d)Q_{i-1}^{s_i/d}
	\]
	the number of monic irreducible polynomials of degree \(s_i\) over
	\(\K_{i-1}\).
	
	By Theorem~\ref{th:relative-root-obstruction}, if two choices
	\[
	G_i(Y),H_i(Y)\in\K_{i-1}[Y]
	\]
	give isotopic semifields at the \(i\)-th step, then their roots satisfy
	\[
	\theta_i'=\kappa\theta_i^\rho
	\]
	for some
	\[
	\kappa\in\K_{i-1}^\times
	\]
	and some
	\[
	\rho\in\operatorname{Aut}(\K_i).
	\]
	Therefore, for a fixed polynomial \(G_i\), the number of polynomials that can
	possibly give an isotopic semifield is at most
	\[
	|\K_{i-1}^\times|\,
	|\operatorname{Aut}(\K_{i-1})|
	=
	(Q_{i-1}-1)|\operatorname{Aut}(\K_{i-1})|.
	\]
	Indeed, the automorphisms of \(\K_i\) fixing \(\K_{i-1}\) only replace a root
	by another conjugate root of the same irreducible polynomial over
	\(\K_{i-1}\), and hence do not produce a new polynomial.
	
	Consequently, after fixing the first \(i-1\) steps, the \(i\)-th lifting step
	produces at least
	\[
	M_i
	:=
	\left\lceil
	\frac{I_{Q_{i-1}}(s_i)}
	{(Q_{i-1}-1)|\operatorname{Aut}(\K_{i-1})|}
	\right\rceil
	\]
	pairwise non-isotopic semifields of order
	\[
	q^{2t_i}
	=
	q^{2t_{i-1}s_i}.
	\]
	
	We now specialize this to the simplest iterative choice.

	\begin{example}
		\label{ex:cubic-iteration-isotopy}
		Assume that
		\[
		s_i=e
		\]
		for every \(i=1,\ldots,k\), where \(e\) is a fixed odd prime. Then
		\[
		t_i=et_{i-1},
		\]
		and hence
		\[
		t_i=te^i.
		\]
		The semifields obtained after \(k\) steps have order
		\[
		N_k=q^{2t_k}=q^{2te^k}.
		\]
		
		At the \(i\)-th step we have
		\[
		Q_{i-1}=q^{t_{i-1}}=q^{te^{i-1}}.
		\]
		Moreover, since \(e\) is prime, the number of monic irreducible
		polynomials of degree \(e\) over \(\K_{i-1}\) is
		\[
		I_{Q_{i-1}}(e)
		=
		\frac{Q_{i-1}^e-Q_{i-1}}{e}.
		\]
		By Theorem~\ref{th:relative-root-obstruction}, two choices
		\(\theta_i,\theta_i'\) can give rise to isotopic semifields only if
		\[
		\theta_i'=\kappa\theta_i^\rho,
		\qquad
		\kappa\in\K_{i-1}^{\times},
		\quad
		\rho\in\Aut(\K_i).
		\]
		Therefore the number of pairwise non-isotopic choices
		at the \(i\)-th step is at least
		\[
		M_i
		=
		\left\lceil
		\frac{Q_{i-1}^e-Q_{i-1}}
		{e\,(Q_{i-1}-1)\,|\operatorname{Aut}(\K_{i-1})|}
		\right\rceil .
		\]
		Write $q=p^f$, where $p$ is prime. Since
\[
|\operatorname{Aut}(\K_{i-1})|=ft e^{i-1},
\]
we have
\[
M_i
=
\left\lceil
\frac{Q_{i-1}+Q_{i-1}^2+\cdots+Q_{i-1}^{e-1}}
     {ft e^i}
\right\rceil
\geq
\frac{q^{te^{i-1}(e-1)}}{ft e^i}.
\]
By Theorem~\ref{thm:equivalence}, liftings of
non-isotopic semifields remain non-isotopic.
Hence the total number of pairwise non-isotopic
semifields obtained after $k$ steps is at least
\[
\prod_{i=1}^k M_i
\geq
\frac{q^{t(e^k-1)}}{(ft)^k e^{k(k+1)/2}}.
\]
		Since
		\[
		N_k=q^{2te^k},
		\]
		we get
		\[
		q^{t(e^k-1)}
		=
		N_k^{\frac{e^k-1}{2e^k}}
		=
		N_k^{\frac12(1-e^{-k})}.
		\]
		Thus, for fixed \(e,t,k\), as \(q\to+\infty\), the construction gives at least
		\[
		N_k^{\frac12(1-e^{-k})+o(1)}
		\]
		pairwise non-isotopic semifields of order \(N_k\), which are
		two-dimensional over their right nucleus. Finally, letting \(k\to+\infty\),
		the exponent
		\[
		\frac12(1-e^{-k})
		\]
		tends to \(1/2\).
	\end{example}

	\section{Newness of the construction}
	\label{sec:comparison-empt-hmo}
	We now compare the semifields in
	Theorem~\ref{th:iterated-explicit-empt-lift} with the known
	constructions. Our aim is to show that the family contains
	semifields which are not isotopic to any previously known
	example. We retain the notation of Section~\ref{sec:explicit-constructions}. The orders of the nuclei provide useful invariants for
distinguishing semifields up to isotopy, although they will
not suffice for all comparisons below. 
    We also recall the Knuth orbit~\cite{knuth1965finite}.
Choose a basis $\{e_1,\ldots,e_n\}$ of $\SSS$ over its center
and write
\[
e_i\star e_j=\sum_{k=1}^n c_{ijk}e_k.
\]
Permuting the three indices of the structure constants $c_{ijk}$
gives presemifields $\SSS^\pi$, with $\pi\in S_3$, called
\emph{Knuth derivatives} of $\SSS$. Their isotopy classes form
the \emph{Knuth orbit}
\[
\mathcal K(\SSS)=\{[\SSS^\pi]:\pi\in S_3\},
\]
which contains at most six classes; here $[\mathbb T]$ denotes
the isotopy class of $\mathbb T$~\cite[Section~2.1]{lavrauw2013tensors}.
Knuth operations preserve the order and the center up to
isomorphism, while permuting the three nuclei up to isomorphism,
computed in semifield representatives of these classes.
For example, the opposite multiplication $x\star^d y=y\star x$
interchanges the left and right nuclei and preserves the middle
nucleus~\cite[Proposition~2.4]{marino2012nuclei}.
Consequently, the multiset
\[
\bigl\{|\NN_l(\SSS)|,|\NN_m(\SSS)|,|\NN_r(\SSS)|\bigr\}
\]
is invariant throughout the Knuth orbit. Distinct such multisets
therefore exclude isotopy even after applying any Knuth operation.

	The comparison with EMPT uses the coordinate sections and
	the unique-line result established in Section~\ref{sec:linear-sets-isotopy}.
	
	\subsection{Comparison with EMPT semifields}
	\label{subsec:comparison-empt}
	\label{subsec:general-lift-noncontainment}

	\begin{proposition}
		\label{prop:explicit-not-empt}
		Let $h\ge1$. No semifield $\SSS_h$ in
		Theorem~\ref{th:iterated-explicit-empt-lift} is isotopic to an
		admissible instance of either construction~\textup{(3)} or
		\textup{(9)} of~\cite{ebert2009infinite}, or to a semifield in their Knuth orbit.
	\end{proposition}

	\begin{proof}
		Put $m=t_h$, $d=t_{h-1}$ and $s=s_h$, so that $m=ds$ and
		$d,s\ge3$ are odd. Set $\E=\F_{q^d}\subseteq \K=\F_{q^m}$ and
		$\LL=\F_{q^{2m}}$. We represent the spreads over their large right
		nuclei, as in Section~\ref{sec:explicit-constructions}.

		We first record an elementary observation. Suppose that
		$\K=\F_q(v)$ and put
		\[
		P_k(v)=\langle1,v,\ldots,v^{k-1}\rangle_{\F_q}.
		\]
		If $k+s\le m$, then $x\E\subseteq P_k(v)$ implies $x=0$.
		Indeed, let $\pi:\K\to\F_q$ extract the coefficient of $v^{m-1}$
		in the power basis. For $c\in \E$ and $0\le j<s$, the inclusion
		$x\E\subseteq P_k(v)$ gives
		\[
		xcv^j\in\langle1,v,\ldots,v^{k+s-2}\rangle_{\F_q},
		\qquad \pi(xcv^j)=0.
		\]
		Since $1,v,\ldots,v^{s-1}$ form an $\E$-basis of $\K$, it follows
		that $\pi(x\K)=0$. This is impossible for $x\ne0$.
		Consequently, if $a\in \LL\setminus \K$, the space
		$P_k(v)\oplus aP_k(v)$ contains no nonzero one-dimensional
		$\E$-subspace: writing $z=x+ay$, the inclusion
		$\E z\subseteq P_k(v)\oplus aP_k(v)$ would imply
		$\E x,\E y\subseteq P_k(v)$.

		For an EMPT semifield of Type~\textup{(3)} isotopic to $\SSS_h$, invariance of
		the order and center forces its defining base field to be
		$\F_q$ and its extension parameter to be $m$; see
		\cite[Proposition~3.1]{ebert2009infinite}.
		Its defining section on the first coordinate line is
		\[
		U_{\rm EMPT}=P\oplus aP,
		\qquad
		P=\langle1,v,\ldots,v^{(m-1)/2}\rangle_{\F_q},
		\]
		where $\K=\F_q(v)$ and $a\notin \K$. With $k=(m+1)/2$, we have
		\[
		k+s=\frac{ds+1}{2}+s\le ds=m.
		\]
		Thus $U_{\rm EMPT}$ contains no nonzero one-dimensional
		$\E$-subspace. In contrast,
		Lemma~\ref{lem:coordinate-sections-iteration} gives
		\[
		\E\theta_h\subseteq\LL_{h-1}\theta_h\subseteq U_h
		\]
		for the first coordinate section of the lifted spread.

		Suppose that these two semifields are isotopic. Their exact
		right nuclei are then both $\K$. By
		Proposition~\ref{prop:isotopy-equivalence-right-nucleus}, the
		isotopy gives a $\K$-semilinear map
		\[
		T:\End_\K(\LL)\longrightarrow\End_\K(\LL),
		\qquad M\longmapsto A M^\rho B,
		\]
		which maps one defining spread subspace exactly onto the other.
		By Corollary~\ref{cor:unique-line-iteration} and
		\cite[Proposition~3.2(iii)]{ebert2009infinite}, the first
		coordinate line is the unique projective line contained in each
		associated linear set. Hence $T$ maps their defining sections on
		these lines onto each other. Every automorphism of $\K$ stabilizes
		its unique subfield $\E$, so this semilinear bijection preserves
		the existence of a nonzero one-dimensional $\E$-subspace.
		This contradicts the preceding paragraph.

		The EMPT construction~\textup{(9)} has order
		$r^{2n}$ and center $\F_r$, where $n$ is even; see the paragraph
		preceding Proposition~3.5 of~\cite{ebert2009infinite}.
		Isotopy would force $r=q$ and $n=m$, contrary to the oddness of
		$m$.

		The nuclear multiset of $\SSS_h$ has a unique large member.
		After placing it on the right, the only Knuth possibilities are
		represented by a spread and its adjoint. For the trace pairing
		$\Tr_{\LL/\K}(xy)$, the adjoint of
		$\varphi_{A,B}(x)=Ax+Bx^{q^m}$ is
		\[
		\varphi_{A,B}^{*}(x)=Ax+B^{q^m}x^{q^m}.
		\]
		Thus the first coordinate section is unchanged, and the same
		obstruction applies.

		We finally check the translation dual. Put $r=(m-1)/2$ and
		$R=\langle1,v,\ldots,v^{r-1}\rangle_{\F_q}$, so that
		$P=\langle1,v,\ldots,v^r\rangle_{\F_q}$.
		If $f$ is the minimal polynomial of $v$ over $\F_q$, then
		orthogonality with respect to $\Tr_{\K/\F_q}(xy)$ gives
		\[
		R^\perp=f'(v)^{-1}P,
		\qquad P^\perp=f'(v)^{-1}R.
		\]
		These identities follow from
		$\Tr_{\K/\F_q}(v^j/f'(v))=0$ for $0\le j\le m-2$
		and comparison of dimensions.
		The coefficient space of Type~\textup{(3)} is
		$C=(P\oplus aP)\times b(R\oplus cR)$.
		The translation dual is obtained by taking its orthogonal
		complement under the polar form of the determinant, followed
		by the trace to $\F_q$; see~\cite{lavrauw2011finite}.
		In these coordinates this form is
		\[
		\beta((A,B),(A',B'))=
		\Tr_{\K/\F_q}\bigl(
		AA'^{q^m}+A^{q^m}A'
		-BB'^{q^m}-B^{q^m}B'\bigr).
		\]
		The two coordinate spaces are orthogonal. Using dual $\K$-bases
		of $\LL$ and the preceding trace identities, their complements
		have the forms
		\[
		(P\oplus aP)^\perp=\mu(R\oplus c'R),
		\qquad
		\bigl(b(R\oplus cR)\bigr)^\perp=\lambda(P\oplus a'P),
		\]
		for some $\lambda,\mu\in \LL^\times$ and $a',c'\notin \K$.
		Swapping the two coordinates and multiplying both by
		$\lambda^{-1}$ gives another admissible spread of Type~\textup{(3)}.
		These operations are spread equivalences: they correspond to
		right composition by $x\mapsto x^{q^m}$ and left multiplication
		by $\lambda^{-1}$. Type~\textup{(3)} is therefore closed under
		translation duality up to isotopy. For Type~\textup{(9)}, both
		power spaces have dimension $n/2$, and the same trace identity
		shows that their orthogonal complements are scalar multiples
		of themselves. The dual-basis argument therefore also proves
		closure of Type~\textup{(9)} under translation duality up to
		isotopy. The preceding exclusions apply to both translation
		duals and their Knuth derivatives.
	\end{proof}
\subsection{Comparison with cyclic and twisted cyclic semifields}
	\label{subsec:comparison-cyclic}

	We now compare the semifields in
Theorem~\ref{th:iterated-explicit-empt-lift} with cyclic
semifields and the twisted cyclic semifields
of~\cite{sheekey2020new}, using the invariance of the
multiset of nuclear orders under Knuth operations.
	Put $m=t_h$. By Theorem~\ref{th:iterated-explicit-empt-lift}, $m>1$ is odd and
	\begin{equation}
		\label{eq:explicit-nuclear-parameters}
		|\SSS_h|=q^{2m},\qquad |Z(\SSS_h)|=q,\qquad
		\bigl(|\NN_l(\SSS_h)|,|\NN_m(\SSS_h)|,|\NN_r(\SSS_h)|\bigr)
		=(q,q,q^m).
	\end{equation}

	\begin{proposition}
		\label{prop:explicit-not-cyclic}
		No semifield $\SSS_h$ in
		Theorem~\ref{th:iterated-explicit-empt-lift} is isotopic to a cyclic
		semifield or to a Knuth derivative of a cyclic semifield.
	\end{proposition}

	\begin{proof}
		By~\cite[Theorem~1]{lavrauw2013semifields}, a cyclic semifield is
		isotopic to a Petit semifield obtained from an irreducible polynomial
		of degree $d$ in $\F_{r^n}[X;\sigma]$, where
		$\Fix(\sigma)=\F_r$. The cases $n=1$ or $d=1$ give fields, whose
		nuclear parameters differ from~\eqref{eq:explicit-nuclear-parameters}.
		For a proper cyclic semifield $\mathbb T$, we have $n,d>1$ and
		\[
		\bigl(|Z(\mathbb T)|,|\NN_l(\mathbb T)|,
		|\NN_m(\mathbb T)|,|\NN_r(\mathbb T)|\bigr)
		=(r,r^n,r^n,r^d)
		\]
		by~\cite[Theorem~9]{lavrauw2013semifields}.
		Hence all three nuclei have order strictly larger than the center.
		This remains true after any Knuth operation, whereas two of the
		nuclei in~\eqref{eq:explicit-nuclear-parameters} have the same order
		as the center.
	\end{proof}

	\begin{proposition}
		\label{prop:explicit-not-twisted-cyclic}
		No semifield $\SSS_h$ in
		Theorem~\ref{th:iterated-explicit-empt-lift} is isotopic to a twisted
		cyclic semifield of~\cite{sheekey2020new} or to a Knuth derivative
		of such a semifield.
	\end{proposition}

	\begin{proof}
		Write $q=p^f$, where $p$ is prime. Suppose that $\SSS_h$ is isotopic
		to a Knuth derivative of a twisted cyclic semifield $\mathbb T$.
		Then $\mathbb T$ has order $q^{2m}$, center of order $q$, and
		multiset of nuclear orders $\{q,q,q^m\}$.

		We use the notation of the skew-polynomial construction, keeping
		its base field distinct from the center of $\mathbb T$.
		Let $r=p^e$, $L=\F_{r^n}$, and let $\sigma\in\Aut(L)$ have order
		$n$ and fixed field $\F_r$. Let $F(Y)\in\F_r[Y]$ be monic
		irreducible of degree $d$, with $F(0)\ne0$. The spread set is the image of
		\[
		\mathcal T=
		\left\{
		\sum_{j=0}^{d-1}a_jX^j+\eta a_0^\rho X^d:
		a_0,\ldots,a_{d-1}\in L
		\right\}
		\]
		in
		\[
		\mathcal A=L[X;\sigma]/(F(X^n))
		\cong M_n(\F_{r^d}),
		\]
		where $\rho\in\Aut(L)$ and the parameters satisfy the hypotheses
		of~\cite{sheekey2020new}. In particular, $|\mathbb T|=r^{nd}$.
		If $\eta=0$, the construction is cyclic, so
		Proposition~\ref{prop:explicit-not-cyclic} applies.
		If $n=1$, the ambient algebra is a field and the spread set has
		the same cardinality as this field, so it defines a field isotope.
		We may therefore assume that $\eta\ne0$ and $n>1$.

		Suppose first that $d\ge3$, and write
		$\rho:x\mapsto x^{p^i}$. The nuclear formulas
		in~\cite[Section~5.1]{sheekey2020new} give
		\[
		\NN_l(\mathbb T)\cong\Fix(\rho),\qquad
		\NN_m(\mathbb T)\cong\Fix(\sigma^d\rho^{-1}),\qquad
		\NN_r(\mathbb T)\cong\F_{r^d},
		\]
		and $Z(\mathbb T)\cong\Fix(\rho|_{\F_r})$.
		Thus $f=\gcd(e,i)$ and $q\le r$. Since $r^d>r\ge q$,
		the right nucleus must have order $q^m$. Comparing this with
		$|\mathbb T|=q^{2m}$ yields $n=2$ and $ed=fm$.
		Consequently $\sigma:x\mapsto x^{p^e}$, and
		\[
		|\NN_l(\mathbb T)|=p^{\gcd(2e,i)},\qquad
		|\NN_m(\mathbb T)|=p^{\gcd(2e,ed-i)}.
		\]
		Write $e=fu$ and $i=fv$. Then $\gcd(u,v)=1$ and $ud=m$,
		so $u$ and $d$ are odd. Since the left nucleus must have order
		$q$, we have $\gcd(2u,v)=1$, and hence $v$ is odd.
		It follows that $ud-v$ is even and
		\[
		\gcd(2u,ud-v)=2,
		\]
		because $\gcd(u,ud-v)=\gcd(u,v)=1$.
		Therefore $|\NN_m(\mathbb T)|=q^2$, contrary to the required
		nuclear orders.

		We treat $d=2$ directly. Write $F(Y)=Y^2+c_1Y+c_0$.
		The element $X$ is invertible in $\mathcal A$, and right
		multiplication of the spread set by $X^{-1}$ gives the equivalent
		spread set
		\[
		\mathcal D=
		\{b+cX^{-1}+\eta c^\rho X:b,c\in L\}.
		\]
		It contains the identity and $L$. A central element is represented
		by an element $z=b+cX^{-1}+\eta c^\rho X$ of $\mathcal D$
		which commutes with every element of $\mathcal D$.
		In $\mathcal A$ we have
		\[
		X^{-1}=-c_0^{-1}(X^{2n-1}+c_1X^{n-1}).
		\]
		For $v\in L$, the coefficient of $X^{2n-1}$ in $zv-vz$ is
		$-c_0^{-1}c(v^{\sigma^{-1}}-v)$.
		Since $n>1$, commutation with every $v\in L$ forces $c=0$.
		Commuting $b$ with $X^{-1}+\eta X\in\mathcal D$ then gives
		$b^\sigma=b$. Thus the center is represented by a subfield of
		$\F_r$, and hence $f\mid e$.

		On the other hand, the scalar field $\F_{r^2}$ of
		$M_n(\F_{r^2})$ centralizes $\mathcal D$ and therefore embeds
		in its right nucleus; see~\cite[Proposition~5]{sheekey2020new}.
		The degree of this right nucleus over the center is consequently
		divisible by $2e/f$, and in particular is even. This contradicts
		the fact that the three nuclear degrees over the center are
		$1,1,m$, all of which are odd.

		Finally, when $d=1$, the construction gives a generalized twisted
		field. In the proper case it is isotopic to the presemifield with multiplication
		\[
		x*y=xy-\eta x^{p^i}y^{p^j},
		\qquad x,y\in\F_{p^{2fm}},
		\]
		where the two Frobenius maps are nonidentity and distinct.
		If either map is the identity, or if the two maps coincide,
		the multiplication is isotopic to that of a field.
		By the nuclear formulas recalled in~\cite[Section~5.1]{sheekey2020new},
		the three nuclear degrees over $\F_p$, up to order, are
		\[
		\gcd(2fm,i),\qquad
		\gcd(2fm,j),\qquad
		\gcd(2fm,j-i).
		\]
		To obtain the multiset $\{f,f,fm\}$, we must have $f\mid i,j$.
		Writing $i=fu$ and $j=fv$, the integers
		\[
		\gcd(2m,u),\qquad
		\gcd(2m,v),\qquad
		\gcd(2m,v-u)
		\]
		would all be odd. Since $2m$ is even, this would require
		$u$, $v$, and $v-u$ all to be odd, which is impossible.
	\end{proof}

\subsection{Proof of newness}
\label{subsec:new-semifields}

To prove that our construction gives a new family, it is enough to
find one member which is not isotopic to any previously known semifield. 

\begin{theorem}
  \label{th:new-semifields-order-2-18}
  \label{prop:general-lift-noncontainment}
  The construction in Theorem~\ref{th:iterated-explicit-empt-lift}
  contains a semifield of order $2^{18}$ which is not isotopic to
  any previously known semifield.
\end{theorem}

\begin{proof}
  Choose an admissible EMPT seed in $\mathcal F_V$ with $q=2$ and
  $t=3$, which exists by~\cite[Theorem~4.5]{ebert2009infinite}.
  Apply Theorem~\ref{th:first-explicit-empt-lift} with $s=3$ and
  any monic irreducible cubic $G(Y)\in\F_8[Y]$. The resulting
  semifield $\SSS_1$ satisfies
  \[
    |\SSS_1|=2^{18},\qquad
    Z(\SSS_1)=\NN_l(\SSS_1)=\NN_m(\SSS_1)\cong\F_2,
    \qquad \NN_r(\SSS_1)\cong\F_{2^9}.
  \]
  Any comparison semifield has nuclear multiset $\{2,2,2^9\}$
  and hence has dimension $2$ over a nucleus.

  Propositions~\ref{prop:explicit-not-cyclic},
  \ref{prop:explicit-not-twisted-cyclic} and~\ref{prop:explicit-not-empt}
  exclude cyclic, twisted cyclic and EMPT semifields, together
  with their Knuth derivatives. The cyclic and twisted cyclic
  comparisons also cover finite fields, generalized twisted fields,
  Sandler and Hughes--Kleinfeld semifields, Knuth Types~II--IV,
  Bierbrauer's construction and its characteristic-two
  generalization, and Dempwolff's construction, including its
  even-characteristic extension; see
  \cite[Sections~4.1 and~5.2--5.4]{koelsch2024unifying}.

  Every HMO semifield has order $r^4$ for a prime power $r$
  \cite[Section~2 and Remark~2.4]{hiramine1987extension}, so its
  order cannot be $2^{18}$. A commutative semifield has equal
  left and right nuclei. Hence a commutative representative in
  the Knuth orbit would have middle nucleus of order $2^9$ and
  left and right nuclei of order $2$. It would be a proper
  commutative semifield of characteristic $2$ and dimension $2$
  over its middle nucleus, contrary to
  \cite[Theorem~3.1]{cohen1982commutative}.

  We check the remaining noncommutative families. A Knuth
  semifield of Type~I over $\K=\F_{2^9}$ has three equal nuclei:
  a nonidentity defining automorphism has odd order, and all
  three nuclei equal its fixed field~\cite[p.~283]{marino2008type}.
  The identity-automorphism case is a field.
  The Johnson--Marino--Polverino--Trombetti constructions and the
  related Marino--Polverino--Trombetti and Marino--Trombetti
  examples have two nuclei of order $r^2$ when the center has
  order $r$~\cite{marino2008type,lavrauw2015solution}; they cannot
  have two nuclei of order $2$.

  Generalized Dickson semifields over $\K=\F_{2^9}$ require a
  nonzero parameter which is not a $(\sigma+1)$st power, where
  $\sigma:x\mapsto x^{2^k}$ is an automorphism of $\K$
  \cite[Section~3]{lavrauw2015solution} and
  \cite[Section~5.1]{koelsch2024unifying}.
  K\"olsch's explicit family in~\cite[Section~6]{koelsch2024unifying}
  requires
  $\alpha\notin(\K^\times)^{\sigma+1}(\K^\times)^{\tau-1}$
  for automorphisms $\sigma,\tau$; see
  \cite[Proposition~4.4]{koelsch2024unifying}.
  Both conditions are impossible, since
  \[
    \gcd(2^k+1,2^9-1)=1
    \qquad(0\le k<9).
  \]
  For a Taniguchi semifield, the nonzero parameter $\eta$ is
  required not to be a $(2^k-1)$st power in $\K$
  \cite[Section~2]{gologlu2024counting}. Consequently,
  $d=\gcd(k,9)>1$. Its multiplication is bilinear over
  $\Fix(\sigma)=\F_{2^d}$, so its unital isotope has center
  containing $\F_{2^d}$, which excludes center $\F_2$.
  The specialization $\eta=0$ lies in the Knuth orbit of Type~I
  \cite[Section~5.6]{koelsch2024unifying}.

  The field-switching construction of~\cite{hou2016switchings}
  is bilinear over its base field, which therefore must be
  $\F_2$. In this case, Theorem~4.1 of that paper forces the
  switching map to be scalar, since the dimension is at least
  $3$. The product is then an invertible additive image of the
  field product and defines a field isotope.

  The other characteristic-two rank-two families of order $r^6$
  in the classical list, including the extensions of the
  Huang--Johnson examples, have center containing $\F_r$
  \cite[Section~3.3.3]{lunardon2022translation}. Order $2^{18}$
  would force $r=8$. The remaining sporadic orders in that list
  differ from $2^{18}$, and the constructions confined to odd
  characteristic do not occur here. 
  The new family of Lobillo, Santonastaso and Sheekey in
  \cite[Theorem~6.4]{lobillo2025quotients} has odd characteristic,
  and their other finite-field specialization is twisted cyclic.
  The semifields in~\cite[Corollary~5.2]{nagy2026semifields}
  have three equal nuclei, while those classified
  in~\cite{gilchrist2026invariants} have order $256$.

  The explicit finite instances in
  \cite[Section~7]{pumpluen2026framework} have order $81$, order
  $r^4$, or order $r^9$ with multiplication bilinear over $\F_r$.
  Only the last order can equal $2^{18}$, but then $r=4$ and
  the center contains $\F_4$.
  Finally, the sufficient norm condition in
  \cite[Theorem~1]{pumpluen2025albert} admits no nonzero twisting
  parameter over $\F_2$, since both sides of the required
  inequality are $1$. A larger base field is contained in the
  center of the normalized semifield.

  These numerical obstructions persist under Knuth operations
  and translation duality. In particular, two nuclei of order
  $2$ force the center to have order $2$, so the central-subfield
  comparisons are also nuclear obstructions. The EMPT translation
  duals are covered by Proposition~\ref{prop:explicit-not-empt}.
  This proves the result.
\end{proof}

\section*{Declaration of AI use}
The majority of this paper was developed and written without the use of generative AI. Generative AI (ChatGPT 6) was only used for investigation purposes in the chapter studying newness (Section 7 of the paper). Outside of that, we used AI only for finding typos,  clearing up notational issues, and minor editorial assistance. 

	\section*{Acknowledgments}
	Paolo Santonastaso is partially supported by the Italian National Group for Algebraic and Geometric Structures and their Applications (GNSAGA- INdAM). Yue Zhou is partially supported by the National Natural Science Foundation of China (No.\ 12371337) and the Natural Science Foundation of Hunan Province (No.\ 2023RC1003).
	
	Paolo Santonastaso and Yue Zhou are deeply grateful to the Kaiyuan International Mathematical Sciences Institute (KIMSI) in Changsha, China, for hosting their visit and providing office space and facilities, where part of this work was completed.

	\bibliographystyle{abbrv}
	\bibliography{biblio}

\end{document}